\documentclass[12pt]{amsart}
\usepackage{amsmath,amssymb,amsfonts,amscd,verbatim,graphicx}
\usepackage{xypic}
\newtheorem{thm}{Theorem}
\newtheorem{lem}[thm]{Lemma}
\newtheorem{conjecture}{Conjecture}

\newtheorem{lemma}[subsection]{Lemma}
\newtheorem{conj}[subsection]{Conjecture}
\newtheorem{proposition}[subsection]{Proposition}

\theoremstyle{definition}
\newtheorem{definition}[subsection]{Definition}

\theoremstyle{remark}

\numberwithin{equation}{section}

\theoremstyle{definition}

\newcommand{\co}{\colon\thinspace}

\begin{document}

\title{A counterexample to the strong version of Freedman's conjecture}
\author{Vyacheslav S. Krushkal}
\address{Department of Mathematics, University of Virginia, Charlottesville, VA 22904}
\email{krushkal\char 64 virginia.edu}

\thanks{Research supported in part by NSF grant DMS-0605280}

\begin{abstract}
A long-standing conjecture due to Michael Freedman asserts that the
4-dimensional topological surgery conjecture fails for non-abelian
free groups, or equivalently that a family of canonical examples of
links (the generalized Borromean rings) are not $A-B$ slice. A
stronger version of the conjecture, that the Borromean rings are not
even weakly $A-B$ slice, where one drops the equivariant aspect of
the problem, has been the main focus in search for an obstruction to
surgery. We show that the Borromean rings, and more generally all
links with trivial linking numbers, are in fact weakly $A-B$ slice.
This result shows the lack of a non-abelian extension of Alexander
duality in dimension $4$, and of an analogue of Milnor's theory of
link homotopy for general decompositions of the $4$-ball.
\end{abstract}

\maketitle

\section{Introduction} \label{introduction}

Surgery and the s-cobordism conjecture, central ingredients of the
geometric classification theory of topological $4-$manifolds, were
established in the simply-connected case and more generally for
elementary amenable groups by Freedman \cite{F0}, \cite{FQ}. Their
validity has been extended to the groups of subexponential growth
\cite{FT}, \cite{KQ}. A long-standing conjecture of Freedman
\cite{F1} asserts that surgery fails in general, in particular for
free fundamental groups. This is the central open question, since
surgery for free groups would imply the general case, cf \cite{FQ}.

There is a reformulation of surgery in terms of the slicing problem
for a special collection of links, the untwisted Whitehead doubles
of the Borromean rings and of a certain family of
their generalizations, see figure 2. (We work in the topological category,
and a link in $S^3=\partial D^4$ is {\em slice} if its components bound disjoint
embedded locally flat disks in $D^4$.) An ``undoubling'' construction \cite{F2}
allows one to work with a more robust link, the Borromean rings, but
the slicing condition is replaced in this formulation by a more
general {\em A--B slice problem}. Freedman's conjecture pinpoints
the failure of surgery in a specific example and states that the
Borromean rings are not $A-B$ slice. This approach to surgery has
been particularly attractive since it is amenable to the tools of
link-homotopy theory and nilpotent invariants of links, and partial
obstructions are known in restricted cases, cf \cite{FL}, \cite{K1},
\cite{K2}. At the same time it is an equivalent reformulation of the
surgery conjecture, and if surgery holds there must exist specific
$A-B$ decompositions solving the problem.

The $A-B$ slice conjecture is a problem at the intersection of
$4-$manifold topology and Milnor's theory of link homotopy \cite{M}.
It concerns codimension zero decompositions of the $4-$ball.
Here a {\em decomposition} of $D^4$, $D^4=A\cup B$, is an extension
of the standard genus one Heegaard decomposition of $\partial
D^4=S^3$. Each part $A, B$ of a decomposition has an attaching
circle (a distinguished curve in the boundary:
${\alpha}\subset\partial A, {\beta}\subset\partial B$) which is the
core of the solid torus forming the Heegaard decomposition of
$\partial D^4$. The two curves ${\alpha}, {\beta}$ form the Hopf
link in $S^3$.

\begin{figure}[ht]
\vspace{.5cm}
\includegraphics[width=3.8cm]{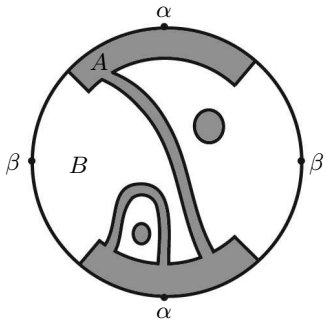}
{\scriptsize
    \put(-59,109){${\alpha}$}
    \put(-59,-5){${\alpha}$}
    \put(-84,89){$A$}
    \put(-116,51){${\beta}$}
    \put(-1,51){${\beta}$}
    \put(-92,50){$B$}}
\vspace{.45cm} \caption{A $2-$dimensional example of a decomposition $(A,{\alpha}),\, (B, {\beta})$:
$D^2=A\cup B$, $A$ is shaded; $({\alpha},\, {\beta})$ are linked $0-$spheres in $\partial D^2$.}
\end{figure}

Figure 1 is a schematic illustration of a decomposition: an example
drawn in two dimensions. While the topology of decompositions in dimension
$2$ is quite simple, they illustrate important basic properties. In this dimension the
attaching regions ${\alpha},\,{\beta}$ are $0-$spheres, and
$({\alpha},\, {\beta})$ form a ``Hopf link'' (two linked $0-$spheres) in $\partial D^2$.
Alexander duality implies that exactly one of the two possibilities holds:
either ${\alpha}$ vanishes as a rational homology class in $A$, or ${\beta}$ does in $B$.
In dimension $2$, this means that either $\alpha$ bounds an arc in $A$, as in the example in figure 1,
or $\beta$ bounds an arc in $B$. (See figure 12 in section \ref{proof} for additional
examples in $2$ dimensions.)

Algebraic and geometric properties of the two parts $A,B$ of a
decomposition of $D^4$ are tightly correlated.
The geometric implication of Alexander
duality in dimension $4$ is that either (an integer multiple of) ${\alpha}$ bounds an
orientable {\em surface} in $A$ or a multiple of $\beta$ bounds a surface in $B$.
Alexander duality does not hold for homotopy groups,
and this difference between being trivial homologically (bounding a surface) as opposed to
{\em homotopically} (bounding a disk) is an algebraic reason for the complexity of
decompositions of $D^4$.

A
geometric refinement of Alexander duality is given by handle structures:
under a mild condition on the handle decompositions which can be assumed without
loss of generality, there is a one-to-one correspondence between
$1-$handles of each side and $2-$handles of its complement. In
general the interplay between the topologies of the two sides is
rather subtle. Decompositions
of $D^4$ are considered in more detail in sections \ref{ABslice} and \ref{decompositions}
of this paper.

We now turn to the main subject of the paper, the $A-B$ slice reformulation
of the surgery conjecture.
An $n-$component link $L$ in $S^3$ is {\em $A-B$ slice} if there
exist $n$ decompositions $(A_i, B_i)$ of $D^4$ and disjoint
embeddings of all $2n$ manifolds $A_1,B_1,\ldots, A_n, B_n$ into
$D^4$ so that the attaching curves ${\alpha}_1,\ldots,{\alpha}_n$
form the link $L$ and the curves ${\beta}_1,\ldots, {\beta}_n$ form
an untwisted parallel copy of $L$. Moreover, the re-embeddings of
$A_i, B_i$ are required to be {\em standard} -- topologically
equivalent to the ones coming from the original decompositions of
$D^4$. The connection of the $A-B$ slice problem for the Borromean
rings to the surgery conjecture is provided by considering the
universal cover of a hypothetical solution to a canonical surgery
problem \cite{F2}, \cite{F3}. The action of the free group by
covering transformations is precisely encoded by the fact that the
re-embeddings of $A_i, B_i$ are standard. A formal definition and a
more detailed discussion of the $A-B$ slice
problem are given in section \ref{ABslice}. The following is the
statement  of Freedman's conjecture \cite{F1}, \cite{F3} concerning
the failure of surgery.

\begin{figure}[ht]
\vspace{.5cm}
\includegraphics[width=3.5cm]{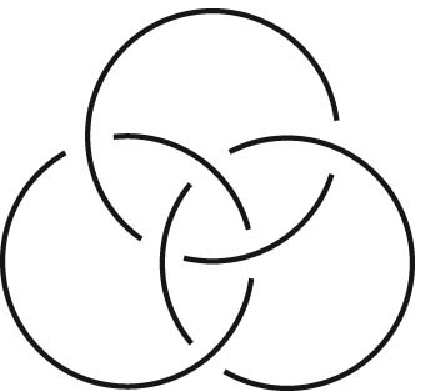} \hspace{2.3cm} \includegraphics[width=3.7cm]{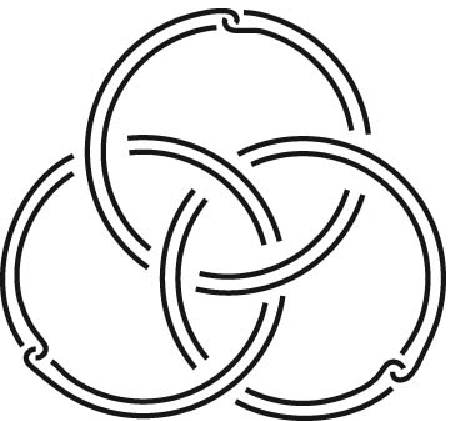}
\vspace{.45cm} \caption{The Borromean rings and their untwisted
Whitehead double.}
\end{figure}

\smallskip

\begin{conjecture} \label{conjecture} \sl The untwisted Whitehead double
of the Borromean rings (figure 2) is not a freely slice link. Equivalently, the
Borromean rings are not $A-B$ slice.
\end{conjecture}

\smallskip

Here a link is {\em freely} slice if it is slice, and in addition
the fundamental group of the slice complement in the $4-$ball is
freely generated by meridians to the components of the link. An
affirmative solution to this conjecture would exhibit the failure of
surgery, since surgery predicts the existence of the free-slice
complement of the link above.

A stronger version of Freedman's conjecture, that the Borromean
rings are not even weakly $A-B$ slice, has been the main focus in
search for an obstruction to surgery. Here a link $L$ is {\em
weakly} $A-B$ slice if the re-embeddings of $A_i, B_i$ are required
to be disjoint but not necessarily standard in the definition above.
To understand the context of this conjecture, consider the simplest
example of a decomposition $D^4=A\cup B$ where $(A,{\alpha})$ is the
$2-$handle $(D^2\times D^2,\partial D^2\times\{0\})$ and $B$ is just the
collar on its attaching curve $\beta$.
This decomposition is trivial in the sense that all topology is
contained in one side, $A$. It is easy to see that a link $L$ is
weakly $A-B$ slice with this particular choice of a decomposition if
an only if $L$ is {\em slice}. The Borromean rings is not
a slice link (cf \cite{M}), so it is not weakly $A-B$ slice with the
trivial decomposition.
However to find an
obstruction to surgery, one needs to find an obstruction for the Borromean
rings to be weakly $A-B$ slice  for all possible decompositions.

Freedman's program in the $A-B$ slice approach to surgery could be
roughly summarized as follows. First consider {\em model}
decompositions, defined using Alexander duality and introduced in
\cite{FL} (see also section \ref{decompositions}). The main step is
then to show that any decomposition is algebraically approximated,
in some sense, by the models -- in this case a suitable algebraic
analogue of the partial obstruction for model decompositions should
give rise to an obstruction to surgery. The first step, formulating
an obstruction for model decompositions, was carried out in
\cite{K2}, \cite{K3}. We now state the main result of this paper
which shows that the second step is substantially more subtle than
previously thought, involving not just the submanifolds but also
their embedding information.

\begin{thm} \label{theorem} \sl Let $L$ be the Borromean rings
or more generally any link is $S^3$ with
trivial linking numbers. Then $L$ is weakly $A-B$ slice.
\end{thm}

The linking numbers provide an obstruction to being weakly $A-B$ slice
(see section \ref{Alexander}), so in fact Theorem \ref{theorem} asserts that a
link is weakly $A-B$ slice if and only if it has trivial linking numbers.

To formulate the main ingredient in the proof of this result in the
geometric context of link homotopy, it is convenient to introduce
the notion of a robust $4-$manifold. Recall that a link $L$ in $S^3$ is
{\em homotopically trivial} if its components bound disjoint
maps of disks in $D^4$. $L$ is called homotopically essential otherwise. (The Borromean
rings is a homotopically essential link \cite{M} with trivial linking numbers.)
Let $(M, {\gamma})$ be a pair
($4-$manifold, attaching curve in $\partial M$). The pair
$(M,{\gamma})$ is {\em robust} if whenever several copies $(M_i,
{\gamma}_i)$ are properly disjointly embedded in $(D^4, S^3)$, the
link formed by the curves $\{ {\gamma}_i\}$  in $S^3$ is
homotopically trivial. The following question relates this notion to
the $A-B$ slice problem: {\sl Given a decomposition $(A, {\alpha}),
(B,{\beta})$ of $ D^4$, is one of the two pairs $(A, {\alpha}),
(B,{\beta})$ necessarily robust?} The answer has been affirmative
for all previously known examples, including the model
decompositions \cite{K2}, \cite{K3}. In contrast, we prove

\begin{lem} \label{theorem 2} \sl There exist decompositions $D^4=A\cup B$
where neither of the two sides $A$, $B$ is robust.
\end{lem}

This result suggests an intriguing possibility that there are
$4-$manifolds which are not robust, but which admit {\em robust
embeddings} into $D^4$. (The definition of a robust embedding $e\co
(M,{\gamma})\hookrightarrow (D^4,S^3)$ is analogous to the
definition of a robust pair above, with the additional requirement
that each of the embeddings $(M_i, {\gamma}_i)\subset (D^4,S^3)$ is
equivalent to $e$.)
 Then the question relevant for the surgery
conjecture is: {\sl given a decomposition $D^4=A\cup B$, is one of
the given {\em embeddings} $A\hookrightarrow D^4$, $B\hookrightarrow
D^4$ necessarily robust?}

Theorem \ref{theorem} has a consequence in the context of {\em
topological arbiters}, introduced in \cite{FK}. Roughly speaking, it
points out a substantial difference in the structure of the
invariants of submanifolds of $D^4$, depending on whether they are
endowed with a specific embedding or not. We refer the reader to
that paper for the details on this application.

Section \ref{ABslice} reviews the background material on surgery and
the $A-B$ slice problem. The $A-B$ slice problem for two-component
links is considered in section \ref{Alexander}; it is shown that
Alexander duality provides an obstruction for links with non-trivial
linking numbers. The proof of theorem \ref{theorem} starts in section
\ref{decompositions} with a construction of the relevant
decompositions of $D^4$. The final section completes the
proof of the theorem.

{\bf Acknowledgements.} This paper concerns the program on the
surgery conjecture developed by Michael Freedman. I would like to
thank him for sharing his insight into the subject on numerous
occasions.

I would like to thank the referee for the comments on the earlier version of this paper.

\section{$4-$dimensional surgery and the the $A-B$ slice problem}
\label{ABslice}

The surgery conjecture asserts that given a $4-$dimensional
Poincar\'{e} pair $(X,N)$, the sequence
$${\mathcal S}^h_{\rm TOP}(X,N)\longrightarrow {\mathcal
N}_{\rm TOP}(X,N)\longrightarrow  L^h_4({\pi}_1 X)$$

is exact (cf [FQ], Chapter 11). This result, as well as the
$5-$dimensional topological s-cobordism theorem, is known to hold
for a class of {\em good} fundamental groups. The simply-connected
case followed from Freedman's disk embedding theorem \cite{F0}
allowing one to represent hyperbolic pairs in ${\pi}_2(M^4)$ by
embedded spheres. Currently the class of good groups is known to
include the groups of subexponential growth \cite{FT}, \cite{KQ} and
it is closed under extensions and direct limits. There is a specific
conjecture for the failure of surgery for free groups \cite{F1}:

\begin{conj} \label{conjecture2}
{\sl There does not exist a topological $4-$manifold $M$, homotopy
equivalent to ${\vee}^3 S^1$ and with $\partial M$ homeomorphic to
${\mathcal S}^0(Wh(Bor))$, the zero-framed surgery on the Whitehead
double of the Borromean rings.}
\end{conj}

This statement is seen to be equivalent to Conjecture
\ref{conjecture} in the introduction by considering the complement
in $D^4$ of the slices for $Wh(Bor)$. This is one of a collection of
canonical surgery problems with free fundamental groups, and solving
them is equivalent to the surgery theorem without restrictions on
the fundamental group. The $A-B$ slice problem, introduced in
\cite{F2}, is a reformulation of the surgery conjecture, and it may
be roughly summarized as follows. Assuming on the contrary that the
manifold $M$ in the conjecture above exists, consider its universal
cover $\widetilde M$. It is shown in \cite{F2} that the end point
compactification of $\widetilde M$ is homeomorphic to the $4-$ball.
The group of covering transformations (the free group on three
generators) acts on $D^4$ with a prescribed action on the boundary,
and roughly speaking the $A-B$ slice problem is a program for
finding an obstruction to the existence of such actions. To state a
precise definition, consider decompositions of the $4-$ball:

\smallskip

\begin{definition} \label{decomposition}
A {\em decomposition} of $D^4$ is a pair of compact
codimension zero submanifolds with boundary $A,B\subset D^4$,
satisfying conditions $(1)-(3)$ below. Denote $$\partial^{+}
A=\partial A\cap
\partial D^4, \; \; \partial^{+} B=\partial B\cap \partial D^4,\; \;
\partial A=\partial^{+} A\cup {\partial}^{-}A, \; \; \partial
B=\partial^{+} B\cup {\partial}^{-}B.$$
(1) $A\cup B=D^4$,\\
(2) $A\cap B=\partial^{-}A=\partial^{-}B,$ \\
(3) $S^3=\partial^{+}A\cup \partial^{+}B$ is the standard genus $1$
Heegaard decomposition of $S^3$.
\end{definition}

Recall the definition of an $A-B$ slice link \cite{F3}, \cite{FL}:

\smallskip

\begin{definition} \label{A-B slice}
Given an $n-$component link $L=(l_1,\ldots,l_n)\subset S^3$, let
$D(L)=(l_1,l'_1,\ldots, l_n,l'_n)$ denote the $2n-$component link
obtained by adding an untwisted parallel copy $L'$ to $L$. The link
$L$ is {\em $A-B$ slice} if there exist decompositions $(A_i, B_i),
i=1,\ldots, n$ of $D^4$ and self-homeomorphisms ${\phi}_i, {\psi}_i$
of $D^4$, $i=1,\ldots,n$ such that all sets in the collection
${\phi}_1 A_1, \ldots, {\phi}_n A_n, {\psi}_1 B_1,\ldots, {\psi}_n
B_n$ are disjoint and satisfy the boundary data:
${\phi}_i({\partial}^{+}A_i)$  is a tubular neighborhood of $l_i$
and ${\psi}_i({\partial}^{+}B_i)$ is a tubular neighborhood of
$l'_i$, for each $i$.
\end{definition}

The surgery conjecture is equivalent to the statement that the
Borromean rings (and a family of their generalizations) are $A-B$
slice. The idea of the proof of one implication is sketched above;
the converse is also true: if the generalized Borromean rings were
$A-B$ slice, consider the complement of the entire collection
${\phi}_i(A_i), {\psi}_i(B_i)$. Gluing the boundary according to the
homeomorphisms, one gets solutions to the canonical surgery problems
(see the proof of theorem 2 in \cite{F2}.)

The restrictions ${\phi}_i|_{A_i}$, ${\psi}_i|_{B_i}$ in the
definition above provide disjoint embeddings into $D^4$ of the
entire collection of $2n$ manifolds $\{A_i, B_i\}$. Moreover, these
re-embeddings are {\em standard}: they are restrictions of
self-homeomorphisms of $D^4$, so in particular the complement
$D^4\smallsetminus {\phi}_i(A_i)$ is homeomorphic to $B_i$, and
$D^4\smallsetminus {\psi}_i(B_i)\cong A_i$. This requirement that
the re-embeddings are standard is removed in the following
definition:

\begin{definition}
A link $L=(l_1,\ldots, l_n)$ in $S^3$ is {\em weakly $A-B$ slice} if
there exist decompositions $((A_1,{\alpha}_1)$, $(B_1,{\beta}_1)),
\ldots, ((A_n,{\alpha}_n), (B_n,{\beta}_n))$ of $D^4$ and disjoint
embeddings of all manifolds $A_i, B_i$ into $D^4$ so that the
attaching curves ${\alpha}_1,\ldots,{\alpha}_n$ form the link $L$
and the curves ${\beta}_1,\ldots, {\beta}_n$ form an untwisted
parallel copy of $L$.
\end{definition}

\section{Abelian versus non-abelian Alexander duality} \label{Alexander}

This section uses Alexander duality to show that the vanishing of the linking numbers
is a necessary condition in theorem \ref{theorem}. Specifically, we prove

\begin{proposition} \label{proposition}
\sl Let $L$ be a link with a non-trivial linking number. Then $L$ is not weakly
$A-B$ slice.
\end{proposition}

{\em Proof.} It suffices to consider $2-$component links, since
any sub-link of a weakly $A-B$ slice link is also weakly $A-B$ slice.
Let  $L=(l_1, l_2)$ with lk$(l_1,l_2)\neq 0$, and consider any two
decompositions $D^4=A_1\cup B_1=A_2\cup B_2$.

Consider the long exact sequences of the pairs $(A_i, \partial^+
A_i), (B_i, \partial^+ B_i)$, where the homology groups are taken
with rational coefficients:
$$ 0\longrightarrow H_2A_i\longrightarrow
H_2(A_i,\partial^+A_i)\longrightarrow
H_1\partial^+A_i\longrightarrow H_1A_i\longrightarrow
H_1(A_i,\partial^+ A_i)\longrightarrow 0$$
$$ 0\longrightarrow H_2B_i\longrightarrow
H_2(B_i,\partial^+B_i)\longrightarrow
H_1\partial^+B_i\longrightarrow H_1B_i\longrightarrow
H_1(B_i,\partial^+ B_i)\longrightarrow 0$$

Recall that $\partial^+A_i, \partial^+B_i$ are solid tori (regular
neighborhoods of the attaching curves ${\alpha}_i, {\beta}_i$.) The
claim is that for each $i$, the attaching curve on exactly one side
vanishes in its first rational homology group. Both of them can't
vanish simultaneously, since the linking number is $1$. Suppose
neither of them vanishes. Then the boundary map in each sequence
above is trivial, and $rk \, H_2(A_i)=rk\, H_2(A_i,\partial^+A_i).$
On the other hand, by Alexander duality $rk \, H_2(A_i)=rk\,
H_1(B_i,\partial^+B_i),\;\, rk \, H_2(A_i,\partial^+A_i)=rk \,
H_1(B_i).$ This is a contradiction, since
$H_1\partial^+B_i\cong{\mathbb Q}$ is in the kernel of
$H_1B_i\longrightarrow H_1(B_i,\partial^+ B_i)$.

Now to show that the link $L=(l_1, l_2)$ is not weakly $A-B$ slice,
set $(C_i,{\gamma}_i)=(A_i,{\alpha}_i)$ if ${\alpha}_i=0\in H_1(A_i;
{\mathbb Q})$ or $(C_i,{\gamma}_i)=(B_i,{\beta}_i)$ otherwise. If
$L$ were weakly $A-B$ slice, there would exist disjoint embeddings
$(C_1,{\gamma}_1)\subset (D^4,S^3)$, $(C_2,{\gamma}_2)\subset
(D^4,S^3)$ so that ${\gamma}_1$ is either $l_1$ or its parallel
copy, and ${\gamma}_2$ is $l_2$ or its parallel copy. Then
$lk({\gamma}_1, {\gamma}_2)\neq 0$, a contradiction. \qed

Proposition \ref{proposition} should be contrasted with theorem \ref{theorem}.
Milnor's link-homotopy invariant of the Borromean rings,
$\overline {\mu}_{123}(Bor)$, equals $1$ \cite{M}.
$\overline {\mu}_{123}$, defined using the quotient
${\pi}_1/({\pi}_1)^3$ of the fundamental group by the third term of the
lower central series, is a non-abelian analogue of the linking number of a link.
Our result, theorem \ref{theorem}, shows the lack of a non-abelian extension of Alexander
duality in dimension $4$.

\section{Decompositions of $D^4$} \label{decompositions}

This section starts the proof of theorem \ref{theorem} by
constructing the relevant decompositions of $D^4$. The simplest
decomposition $D^4=A\cup B$ where $A$ is the $2-$handle $D^2\times
D^2$ and $B$ is just the collar on its attaching curve, was
discussed in the introduction. Now consider the genus one surface
$S$ with a single boundary component ${\alpha}$, and set $A_0=S\times D^2$.
Moreover, one has to specify its embedding into
$D^4$ to determine the complementary side, $B$. Consider the
standard embedding (take an embedding of the surface in $S^3$, push
it into the $4-$ball and take a regular neighborhood.) Note that
given any decomposition, by Alexander duality the attaching curve of
exactly one of the two sides  vanishes in it homologically, at least
rationally. Therefore the decomposition under consideration now may
be viewed as the first level of an ``algebraic approximation'' to an
arbitrary decomposition.

\begin{figure}[ht]
\vspace{.5cm}
\includegraphics[width=2.8cm]{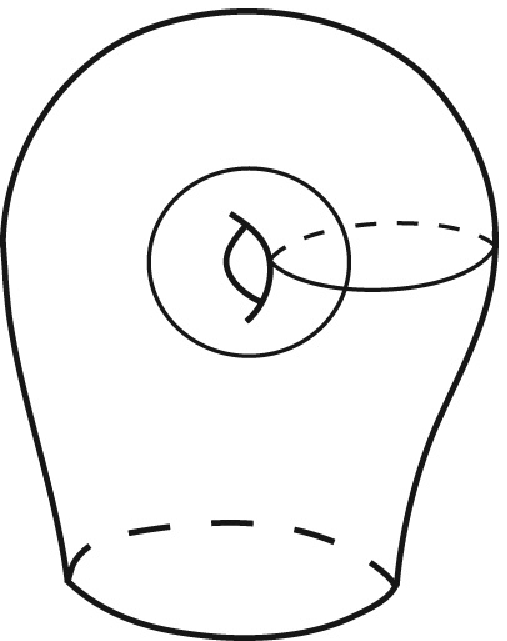} \hspace{3cm} \includegraphics[width=3.8cm]{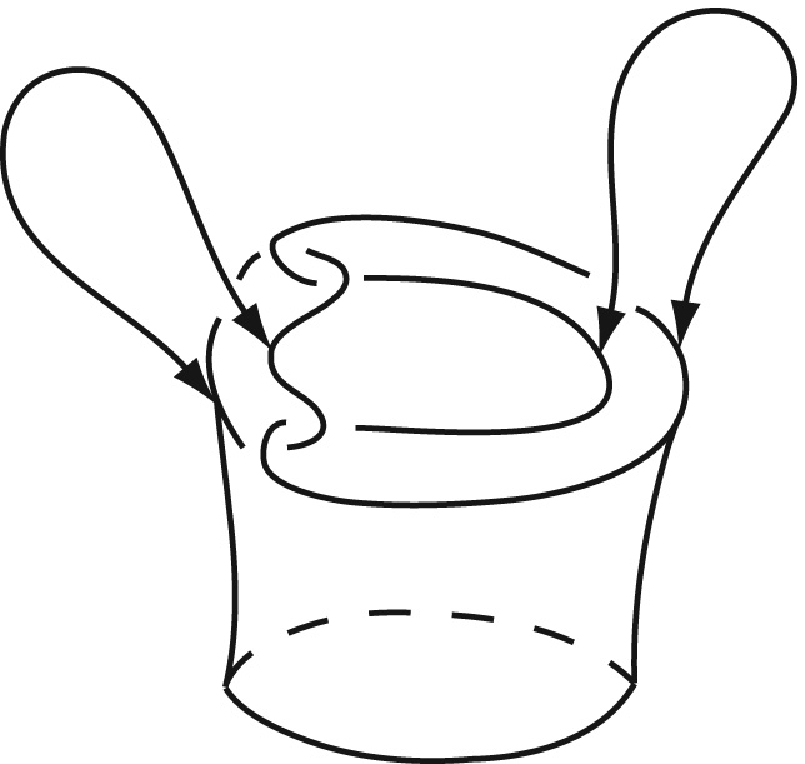}
{\small
    \put(-297,27){$A_0$}
    \put(-264,41){${\alpha}_1$}
    \put(-221,46){${\alpha}_2$}
    \put(-270,-7){${\alpha}$}
    \put(-80,-7){${\beta}$}
    \put(-126,73){$H_1$}
    \put(2,80){$H_2$}
    \put(-8,23){$B_0$}}
\vspace{.45cm} \caption{}
\end{figure}

\begin{proposition} \label{surface complement}
\sl Let $A_0=S\times D^2$, where $S$ is the genus one
surface with a single boundary component $\alpha$. Consider the
standard embedding $(A_0, {\alpha}\times\{ 0\})\subset (D^4,
S^3)$. Then the complement $B_0$ is obtained from the collar
on its attaching curve, $S^1\times D^2\times I$, by attaching a pair
of zero-framed $2-$handles to the Bing double of the core of the
solid torus $S^1\times D^2\times\{1\}$, figure 3.
\end{proposition}

\begin{figure}[ht]
\vspace{.54cm}
\includegraphics[width=4cm]{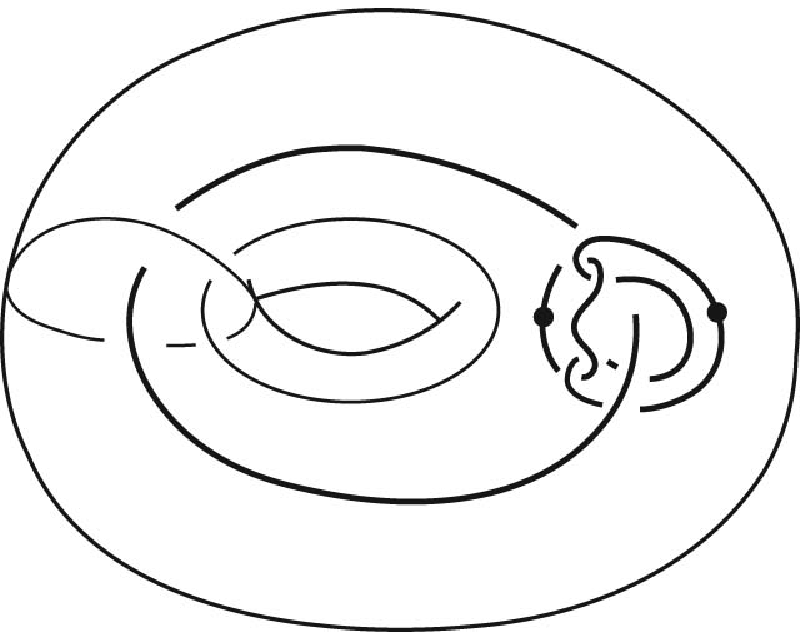} \hspace{3cm} \includegraphics[width=4cm]{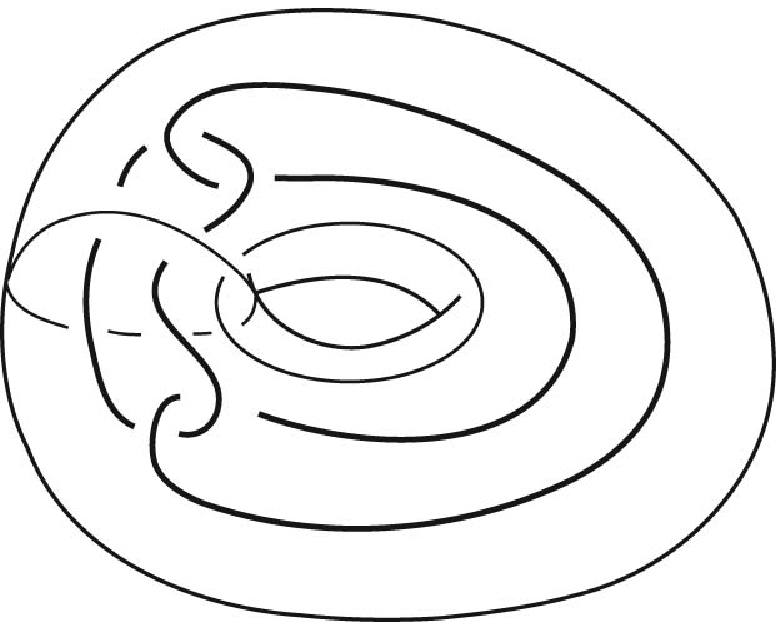}
{\small
    \put(-336,5){$A_0$}
    \put(-128,5){$B_0$}}
{\scriptsize
    \put(-304,20){$0$}
    \put(-260,29){$\alpha$}
    \put(-107,30){$0$}
    \put(-17,22){$0$}
    \put(-46,35){$\beta$}}
\vspace{.48cm} \caption{}
\end{figure}

The proof is a standard exercise in Kirby calculus, see for example
\cite{FL}. A precise description of these $4-$manifolds is given in
terms of Kirby diagrams in figure $4$. Rather than considering
handle diagrams in the $3-$sphere, it is convenient to draw them in
the solid torus, so the $4-$manifolds are obtained from $S^1\times
D^2\times I$ by attaching the $1-$ and $2-$handles as shown in the
diagrams. To make sense of the ``zero framing'' of curves which are
not null-homologous in the solid torus, recall that the solid torus
is embedded into $S^3=\partial D^4$ as the attaching region of a
$4-$manifold, and the $2-$handle framings are defined using this
embedding.

This example illustrates the general principle that (in all examples
considered in this paper) the $1-$handles of each side are in
one-to-one correspondence with the $2-$handles of the complement.
This is true since the embeddings in $D^4$ considered here are all
standard, and in particular each $2-$handle is unknotted in $D^4$.
The statement follows from the fact that $1-$handles may be viewed
as standard $2-$handles removed from a collar, a standard technique
in Kirby calculus (see Chapter 1 in \cite{Ki}.)
Moreover, in each of our examples the attaching
curve ${\alpha}$ on the $A-$side bounds a surface in $A$, so it has
a zero framed $2-$handle attached to the core of the solid torus. On
the $3-$manifold level, the zero surgery on this core transforms the
solid torus corresponding to $A$ into the solid torus corresponding
to $B$. The Kirby diagram for $B$ is obtained by taking the diagram
for $A$, performing the surgery as above, and replacing all zeroes
with dots, and conversely all dots with zeroes. (Note that the
$2-$handles in all our examples are zero-framed.)

Note that a distinguished pair of curves ${\alpha}_1, {\alpha}_2$,
forming a symplectic basis in the surface $S$, is determined as the
meridians (linking circles) to the cores of the $2-$handles $H_1,
H_2$ of $B_0$ in $D^4$. In other words, ${\alpha}_1$,
${\alpha}_2$ are fibers of the circle normal bundles over the cores
of $H_1, H_2$ in $D^4$.

\begin{figure}[ht]
\vspace{.69cm}
\includegraphics[width=3.8cm]{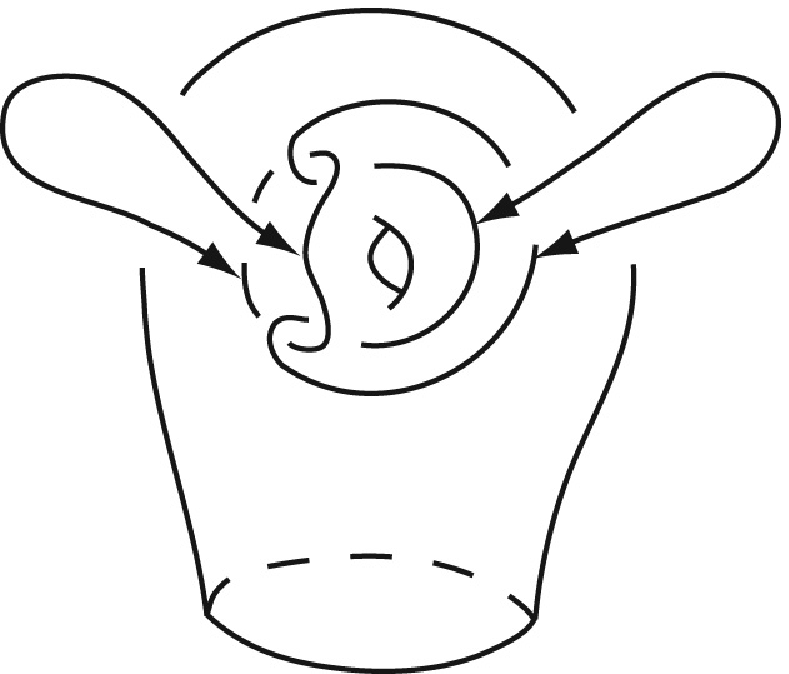}
 \hspace{2.5cm} \includegraphics[width=6cm]{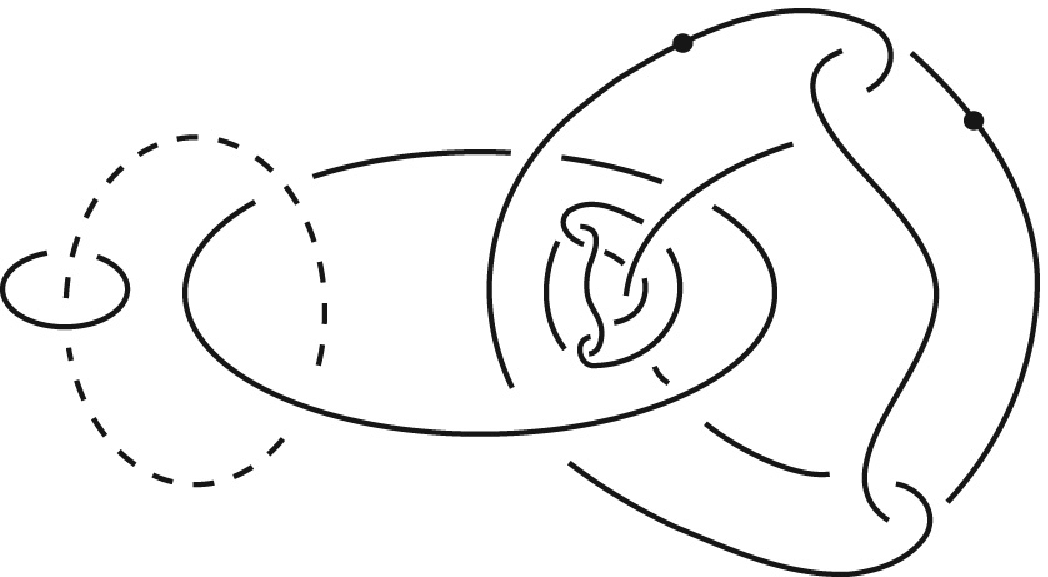}
{\small
    \put(-358,40){$A_1$}
    \put(-340,0){${\alpha}$}
    \put(-182,47){${\alpha}$}}
{\scriptsize
    \put(-99,13){$0$}
    \put(-88,47){$0$}
    \put(-57,47){$0$}}
    \vspace{.5cm} \caption{}
\end{figure}

\begin{figure}[ht]
\vspace{.5cm}
\includegraphics[width=3.8cm]{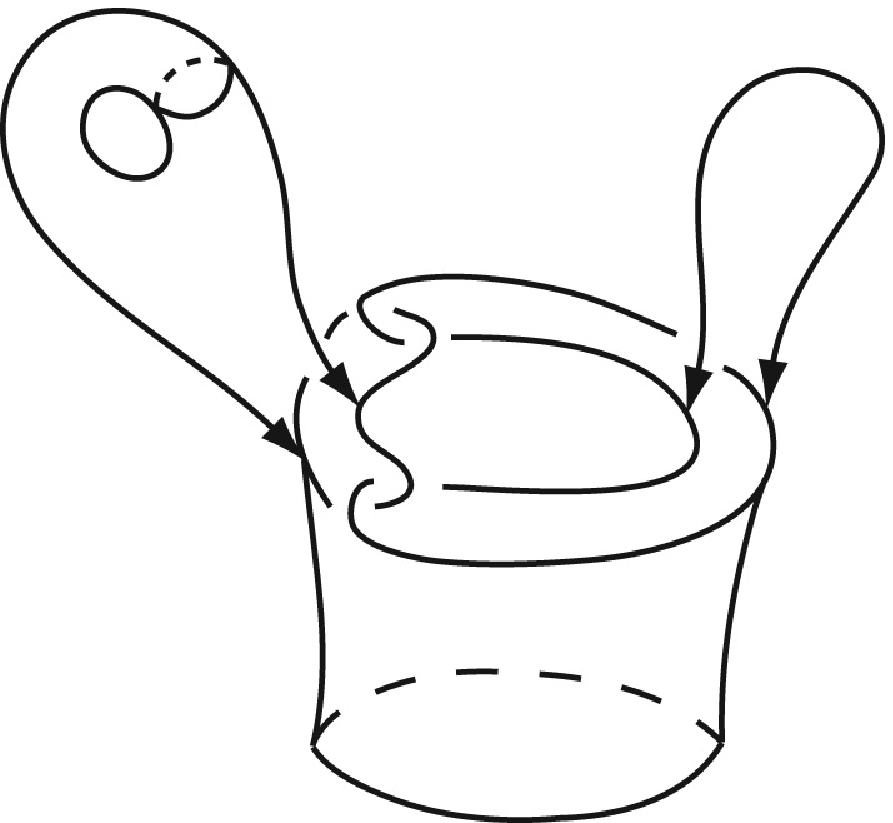} \hspace{2.2cm} \includegraphics[width=6cm]{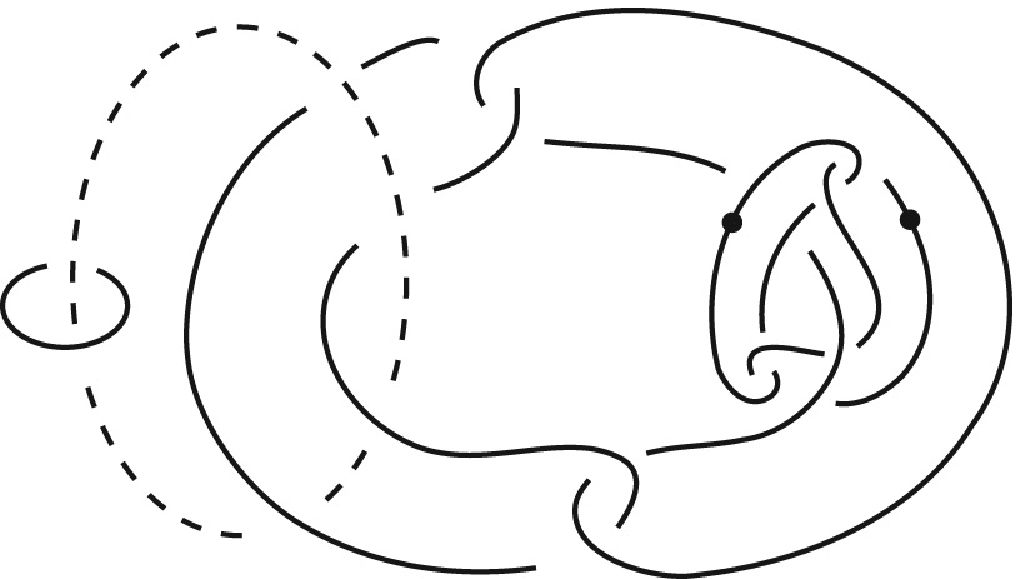}
{\small
    \put(-346,34){$B_1$}
    \put(-330,0){${\beta}$}
    \put(-184,42){${\beta}$}}
{\scriptsize
    \put(-101,-10){$0$}
    \put(-61,-11){$0$}}
    \vspace{.45cm} \caption{}
\end{figure}

An important observation \cite{FL} is that this construction may be
iterated: consider the $2-$handle $H_1$ in place of the original
$4-$ball. The pair of curves (${\alpha}_1$, the attaching circle
${\beta}_1$ of $H_1$) form the Hopf link in the boundary of $H_1$.
In $H_1$ consider the standard genus one surfaces bounded by
${\beta}_1$. As discussed above, its complement is given by two
zero-framed $2-$handles attached to the Bing double of ${\alpha}_1$.
Assembling this data, consider the new decomposition $D^4=A_1\cup
B_1$, figures 5, 6. As above, the diagrams are drawn in solid tori
(complements in $S^3$ of unknotted circles drawn dashed in the
figures.) The handlebodies $A_1, B_1$ are examples of {\em model
decompositions} \cite{FL} obtained by iterated applications of the
construction above. It is shown in \cite{K2}, \cite{K3} that such
model handlebodies are robust, or in other words the Borromean rings
are not weakly $A-B$ slice when restricted to the class of model
decompositions.

\begin{figure}[ht]
\vspace{.5cm}
\includegraphics[width=3.8cm]{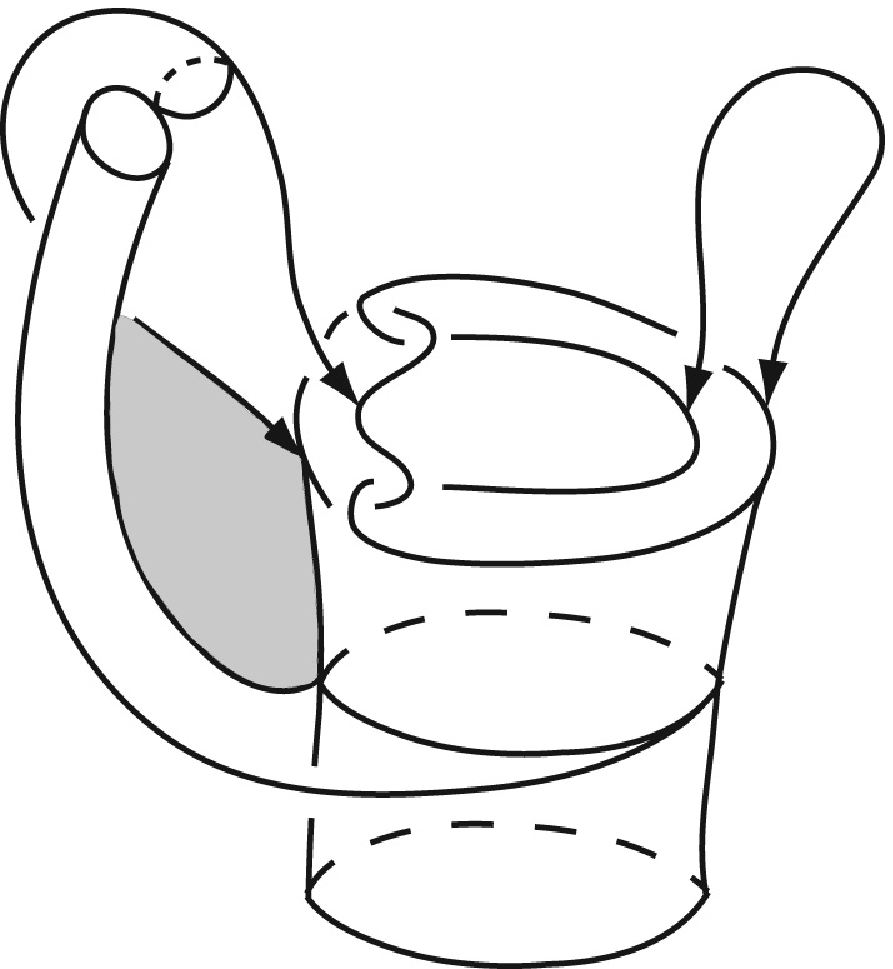} \hspace{2.4cm} \includegraphics[width=6.2cm]{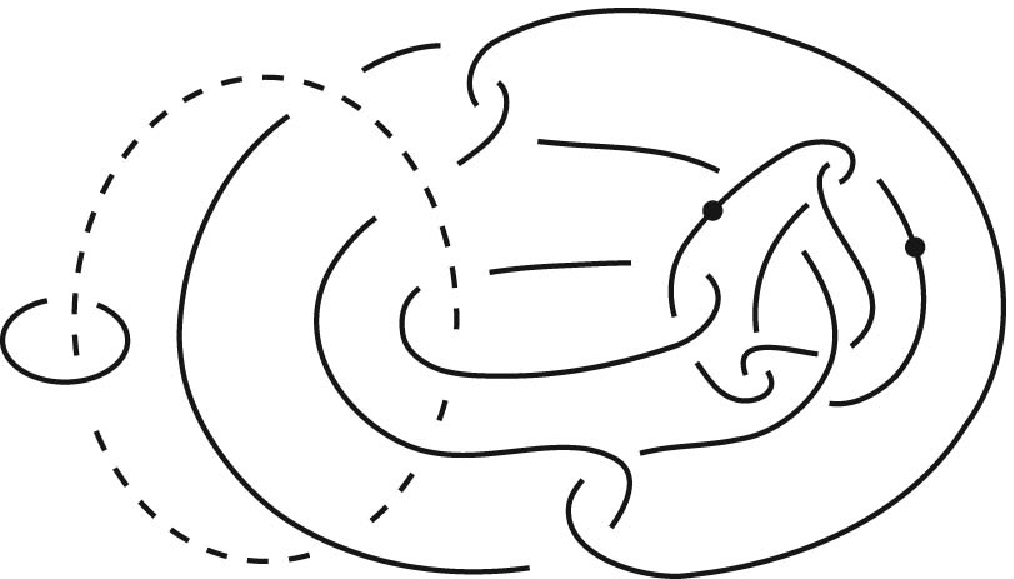}
{\small
    \put(-366,30){$B$}
    \put(-331,-4){${\beta}$}
    \put(-187,35){$\beta$}}
{\scriptsize
    \put(-98,-10){$0$}
    \put(-58,-11){$0$}
    \put(-80,59){$0$}}
    \vspace{.45cm} \caption{}
\end{figure}

We are now in a position to define the decomposition $D^4=A\cup B$
used in the proof of theorem \ref{theorem}.

\begin{definition} \label{definition} Consider
$B=(B_1\, \cup $ zero-framed $2-$handle$)$ attached as shown in the
Kirby diagram in figure 7. The effect of this $2-$handle on the
complement $A=D^4\smallsetminus B$ is shown in figure 8: it adds a
$1-$handle to the diagram of $A_1$. Figure 9 shows a handle diagram
of $A$ after a handle slide. Note that a $(1-, 2-)$ handle pair may
be canceled, the result is given on the left in figure 12. This fact
will be used in the proof of theorem \ref{theorem}.
\end{definition}

\begin{figure}[ht]
\vspace{.5cm}
\includegraphics[width=4.4cm]{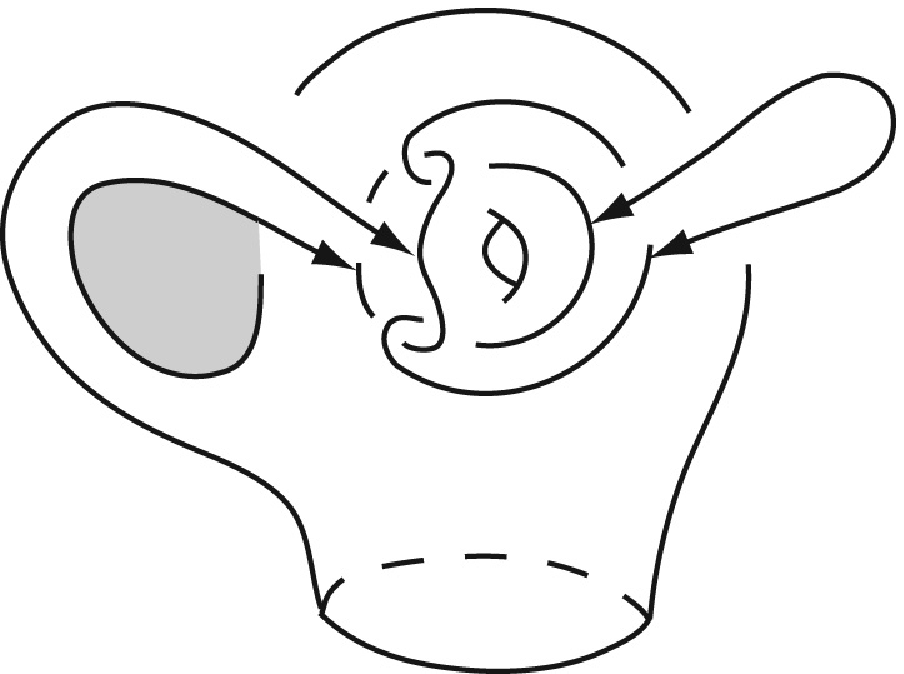} \hspace{2.4cm} \includegraphics[width=6.2cm]{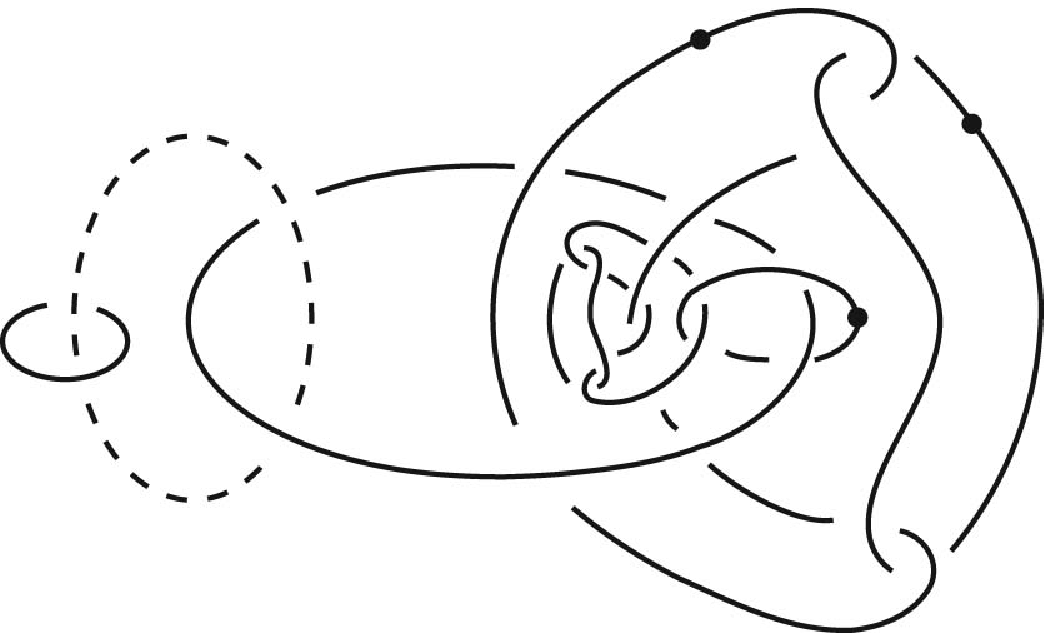}
{\small
    \put(-375,26){$A$}
    \put(-343,0){${\alpha}$}
    \put(-187,46){${\alpha}$}}
{\scriptsize
    \put(-105,17){$0$}
    \put(-90.5,49.5){$0$}
    \put(-61,37.5){$0$}}
    \vspace{.45cm} \caption{}
\end{figure}

\begin{figure}[ht]
\vspace{.5cm}
\includegraphics[width=6.2cm]{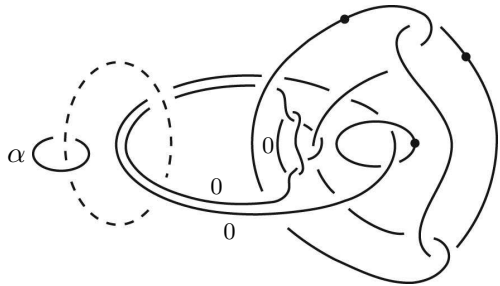}
{\small
    \put(-187,46){${\alpha}$}}
{\scriptsize
    \put(-105,17){$0$}
    \put(-90.5,49.5){$0$}
    \put(-110,34){$0$}}
    \vspace{.45cm} \caption{A handle diagram for $A$ after a handle slide.}
\end{figure}

Imprecisely (up to homotopy, on the level of spines) $B$ may be
viewed as $B_1\cup 2-$cell attached along (the attaching circle
${\beta}$ of $B_1$, followed by a curve representing a generator of
$H_1$ of the second stage surface of $B_1$). This $2-$cell is
schematically shown in the spine picture of $B$ in the first part of
figure 7 as a cylinder connecting the two curves. The shading
indicates that the new generator of ${\pi}_1$ created by adding the
cylinder is filled in with a disk. Similarly, one checks that the
effect of this operation on the $A-$side is that one of the
$2-$handles at the second stage is connected-summed with the first
stage surface, figure 8. (This is seen in the handle diagram by
canceling a $1-, 2-$handle pair, as shown in figure 12.) Again, the shading indicates that
no new generators of ${\pi}_1$ are created. The figures showing the
spines are provided only as a motivation for the construction; a
precise description of $A, B$ is of course given by their handle
diagrams. While the proof of theorem \ref{theorem} below is given in
terms of Kirby diagrams, it can easily be followed at the level of
spines.

\section{Proof of theorem \ref{theorem}: a relative slice problem} \label{proof}

We start this section by recalling the technique which will be
useful in completing the proof of theorem \ref{theorem}, the {\em
relative slice problem}, introduced in \cite{FL}. The setup in our
context is as follows: suppose two codimension zero submanifolds $M,
N$ of $D^4$ are given; each one has an attaching circle
${\gamma}\subset
\partial M$, $\delta\subset\partial N$. The submanifolds are
proper in the sense that one has embeddings of pairs
$(M,{\gamma})\subset (D^4, S^3)$, $(N, {\delta})\subset (D^4,S^3)$,
where each circle ${\gamma}$, ${\delta}$ is unknotted in the
$3-$sphere.

The problem that has to be analyzed is: can $(M,{\gamma})$,
$(N,{\delta})$ be embedded {\em disjointly} into $(D^4, S^3)$ so
that the curves ${\gamma}, {\delta}$ form the Hopf link in the
$3-$sphere? Assume that $M, N$ have handle decompositions, relative
to the attaching regions $S^1\times D^2$, with only $1-$ and
$2-$handles. Let ${\gamma}, {\delta}$ form the Hopf link in
$\partial D^4$, and consider the $4-$ball
$D'=D^4\smallsetminus$(collar on $\partial D^4$). To be
precise, denote the $1-$handles of $M, N$ by ${\mathcal H}_1$,
${\mathcal H}'_1$, and their $2-$handles by ${\mathcal H}_2$,
${\mathcal H}'_2$. As usual, we view the 1-handles of $M,N$ as
standard slices removed from their collars. Denote these slices by
${\mathcal H}^*_1$, ${\mathcal H}'^*_1$. Then $M, N$ embed
disjointly into $D^4$ if and only if there are disjoint embeddings
of the $2-$handles ${\mathcal H}_2\cup {\mathcal H}'_2$, attached to
the collars, in the handlebody $D'\cup{\mathcal H}^*_1\cup{\mathcal
H}'^*_1$.

\begin{figure}[ht]
\vspace{.4cm}
\includegraphics[height=4cm]{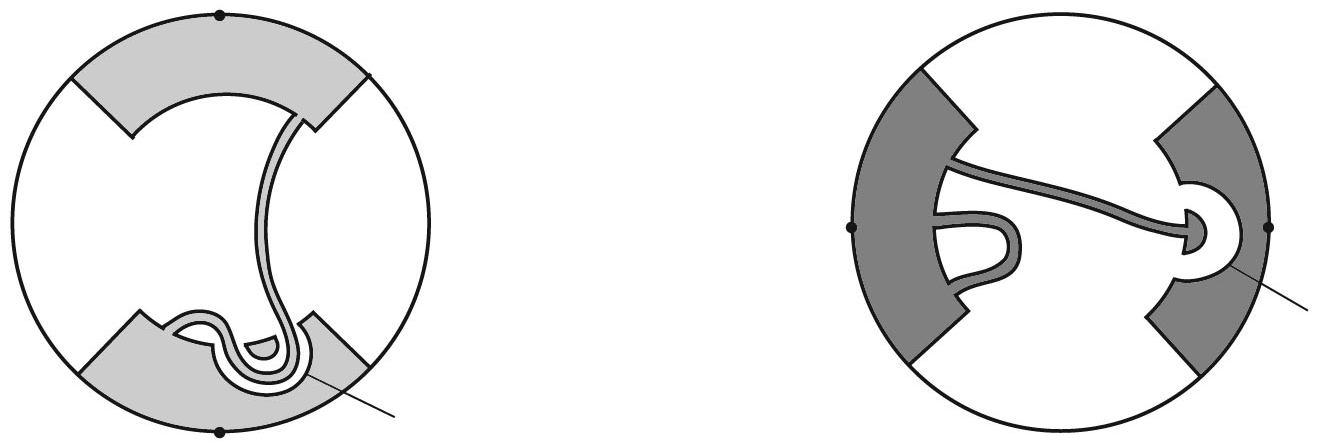}
{\scriptsize
    \put(-284,116){${\gamma}$}
    \put(-284,-6){${\gamma}$}
    \put(-235,3){${\mathcal H}^*_1$}
    \put(-287,52){${\mathcal H}_2$}
    \put(-305,93){$M$}
    \put(-127,53){${\delta}$}
    \put(-9,53){${\delta}$}
    \put(0,30){${\mathcal H}'^*_1$}
    \put(-72,50){${\mathcal H}'_2$}
    \put(-110,37){$N$}}
\vspace{.4cm} \caption{}
\end{figure}

\begin{figure}[ht]
\includegraphics[height=4cm]{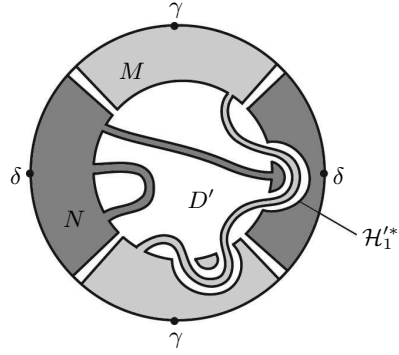}
{\scriptsize
    \put(-72,118){${\gamma}$}
    \put(-72,-7){${\gamma}$}
    \put(-91,93){$M$}
    \put(-132,53){${\delta}$}
    \put(-10,53){${\delta}$}
    \put(1,30){${\mathcal H}'^*_1$}
    \put(-112,37){$N$}
    \put(-65,44){$D'$}}
 \vspace{.4cm} \caption{Disjoint embeddings of $(M,{\gamma}),\, (N, {\delta})$ in figure 10 into $(D^4,S^3)$,
 where ${\gamma}, \, {\delta}$ form a Hopf link in $S^3$.}
\end{figure}

An example of $M,\, N$ drawn in two dimensions is given in figure 10, and a
solution to this relative-slice problem -- disjoint embeddings of
$M,\, N$ in $D^4$ with their attaching circles ${\gamma},\, {\delta}$
forming a Hopf link in $\partial D^4$ -- is shown in figure 11. Note
that the handle ${\mathcal H}_2$ of $M$ in the solution goes over
the ``helping'' handle ${\mathcal H}'^*_1$ attached to $D'$.

Consider the decomposition $D^4=A\cup B$ constructed in definition
\ref{definition}. The proof of theorem \ref{theorem} follows from
lemmas \ref{lemma A} and \ref{lemma B} below.

\smallskip

\begin{lemma} \label{lemma A}
\sl Let $S$ denote the genus one surface with one boundary
component, ${\gamma}=\partial S$. Denote by $S_0$ its
untwisted $4-$dimensional thickening, $S_0=S\times D^2$, and
set ${\gamma}_0={\gamma}\times\{0\}$. Then there exists a
proper embedding $(A, {\alpha})\subset (S_0,
{\gamma}_0)$.
\end{lemma}

{\em Proof}. Kirby diagrams of $A$ are given in figures 8, 9. Observe
that a ($1-$, $2-$handle) pair in the diagram in figure 9 may be
canceled, the result is shown on the left in figure 12.

\begin{figure}[ht]
\vspace{.5cm}
\includegraphics[height=3.5cm]{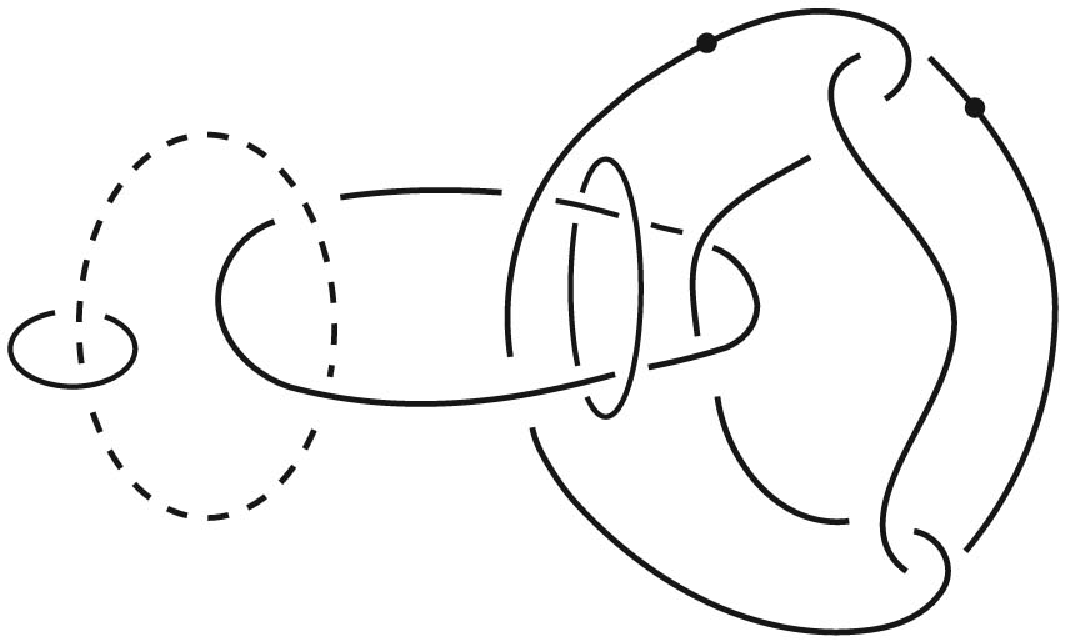} \hspace{1.7cm} \includegraphics[height=3.5cm]{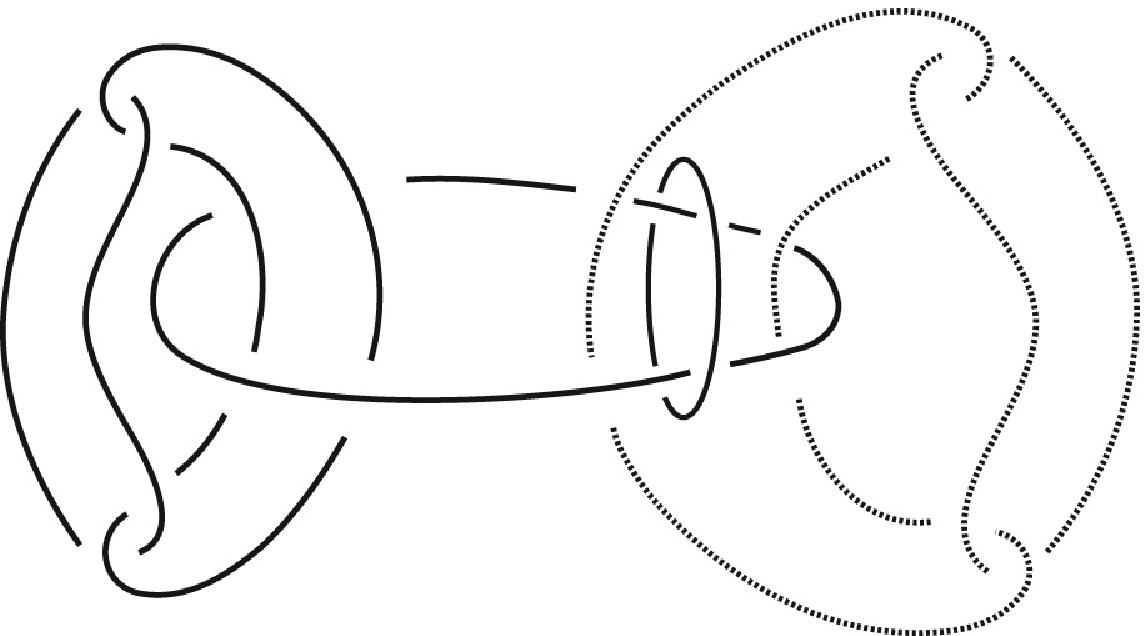}
{\scriptsize
    \put(-308,25){$0$}
    \put(-336,26){$0$}
    \put(-185,73){$l_1$}
    \put(-136,90){$l_2$}
    \put(-108,25){$l_3$}
    \put(-69,25){$l_4$}
    \put(-5,80){$r_2$}
    \put(-71,93){$r_1$}
    \put(-355,2){$A$}
    \put(-400,43){${\alpha}$}}
\vspace{.45cm} \caption{}
\end{figure}

In light of proposition \ref{surface complement}, to prove that $A$
embeds in $S_0$ it suffices to show that $(A, {\alpha})$
embeds in the complement of a standard embedding of two zero-framed
$2-$handles attached to the Bing double of a meridian to $\alpha$ in
$S^3$. This is an instance of the relative-slice problem discussed
above, where $(M,{\gamma})=(A,{\alpha})$ and $N$ is obtained from a
collar on ${\delta}$ by attaching $2-$handles to the Bing double of
the core. (Note $(N,{\delta})$ equals $(B_0,{\beta})$ considered in
sections \ref{decompositions}, see figures 3, 4.)
This relative-slice problem is shown on the right in figure 12. The link
is considered in the $3-$sphere boundary of the $4-$ball $D'$, and
the link $l_1, \ldots, l_4$ has to be sliced in the handlebody
$D'\cup (2-$handles ${\mathcal H}^*_1)$ where the handles are
attached with zero framings
along $r_1, r_2$. Here $l_1,l_2$ are the attaching curves for the
$2-$handles of $N$ and $l_3, l_4$ are the attaching curves for the
$2-$handles of $M$. Note that the slices for $l_1, l_2$ constructed
in the proof are required to be standard in $D^4$, to make sure that
their complement is the thickened surface $S_0$.

A solution to this relative-slice problem is given in figures 13, 14.
The slices are described in terms of the Morse function given by the radial
coordinate in the $4-$ball $D'$. Denote the $3-$sphere at the radius $R$ from
the origin by $S^3_R$, $0<R\leq 1$. The link on the right in figure 12 lies in
$\partial D'=S^3_1$. The link components move by an isotopy
for $1 > R> 3/4$, and at $R=3/4$ the component $l_4$ is
connected-summed with a parallel copy of $r_2$. The result is denoted by $l'_4$,
figure 13. Note that $l'_4$ bounds a disk in $S^3_{3/4}$ in the
complement of all other curves. To make the slice non-degenerate in terms of
the Morse function, let $l'_4$ bound a disk as $R$ decreases from $3/4$ to $1/2$,
while all other
curves move by an isotopy. The link in $S^3_{1/2}$ is shown on the right in figure 13.

\begin{figure}[ht]
\vspace{.5cm}
\includegraphics[width=6.5cm]{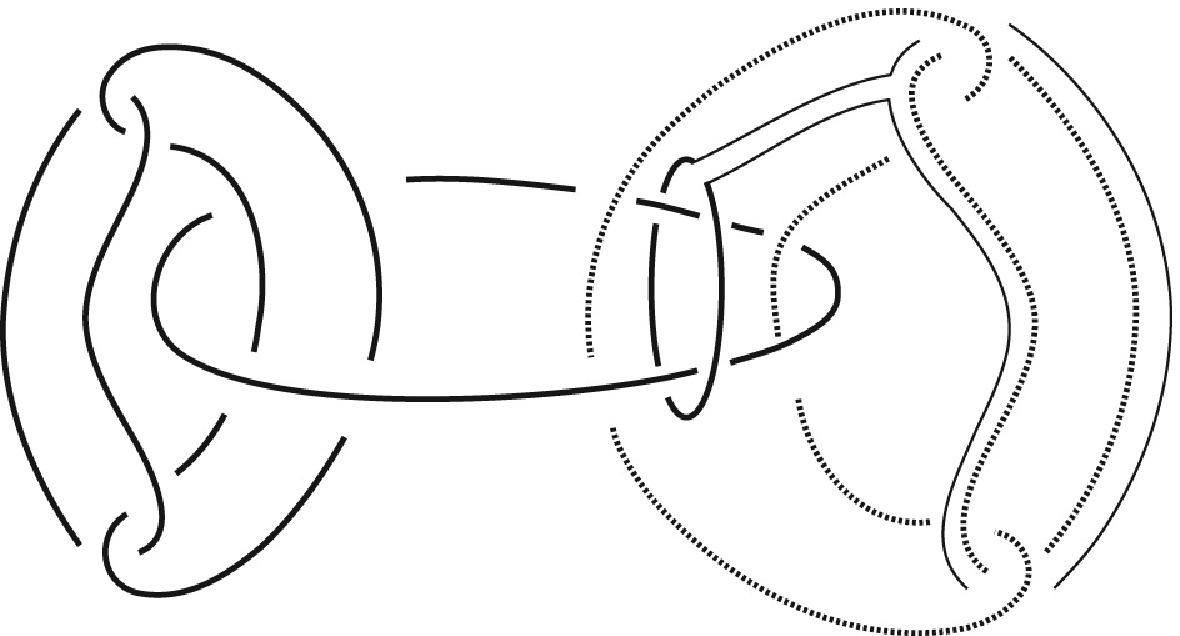} \hspace{1.2cm} \includegraphics[width=6.5cm]{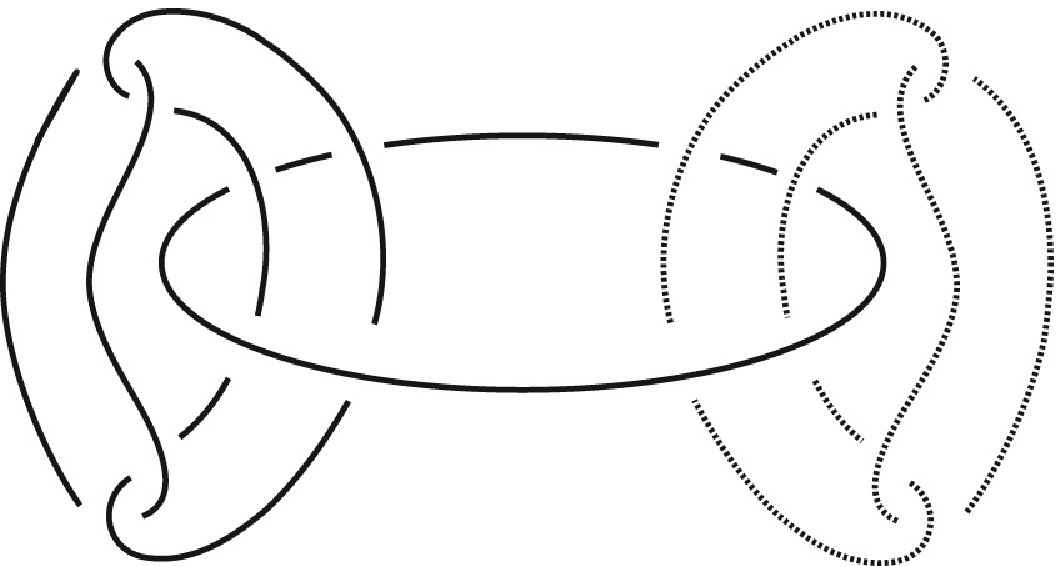}
{\scriptsize
    \put(-189,73){$l_1$}
    \put(-132,90){$l_2$}
    \put(-105,21){$l_3$}
    \put(-299,25){$l'_4$}
    \put(-3,80){$r_2$}
    \put(-233,80){$r_2$}
    \put(-61,93){$r_1$}
    \put(-335,-3){R=3/4}
    \put(-84,-3){R=1/2}}
\vspace{.5cm} \caption{}
\end{figure}

The curves $r_i\subset \partial D'$ bound disjoint embedded disks ${\Delta}_i$:
the cores of the zero-framed $2-$handles ${\mathcal H}^*_1$ attached to $D'$. As the Morse function $R$ changes
from $1$ to $0$, it is important to note that the curves $r_i$ move by an isotopy
and no other curves intersect them. Therefore, $r_1,\, r_2$ at each radius $R_0$ bound
disjoint disks: the disks ${\Delta}_i$ as above, union with the annuli corresponding
to the isotopy of $r_i$ for $1>R>R_0$. Moreover, since the handles ${\mathcal H}^*_1$
attached to $D'$ are zero-framed, untwisted parallel copies of $r_i$ also bound disjoint
embedded disks.

Morse-theoretically the connected sum at $R=3/4$ in figure 13 corresponds
to a saddle point of the slice for $l_4$. This slice is of the form shown
in figure 15 (disregard the labels in that figure, which are used for a later argument.)

To finish the proof of the relative-slice problem, let the link in $S^3_{1/2}$ move
by an isotopy for $1/2>R>1/4$, and at $R=1/4$ the components $l_1,\, l_2$ are connected-summed
with $r_1,\, r_2$ as shown on the left in figure 14. Denote the resulting curves by $l'_1,l'_2$.
The components $l'_1, l'_2, l_3$ form the unlink. This is seen by performing an isotopy
(at $1/4>R>1/8$) to the link on the right in figure 14. Now let all curves bound disks at
$1/8>R>0$. The slices for $l_1,\, l_2$ again have the form shown in figure 15; the slice for
$l_3$ has just a single critical point.

\begin{figure}[ht]
\vspace{.5cm}
\includegraphics [width=5.5cm]{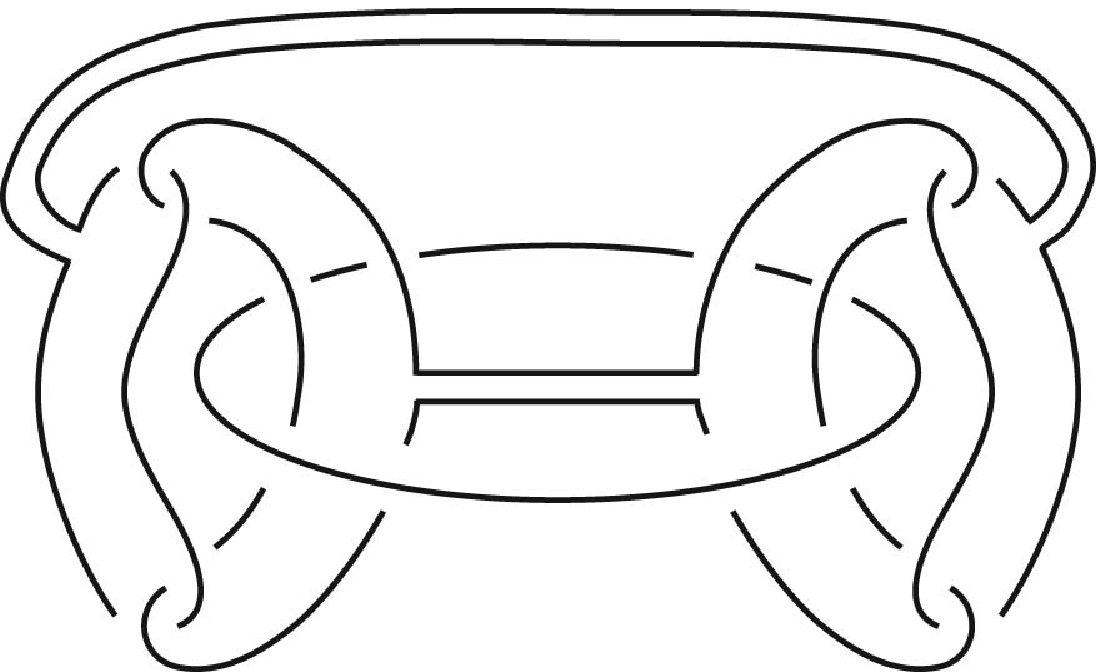} \hspace{2.5cm} \includegraphics[width=4.5cm]{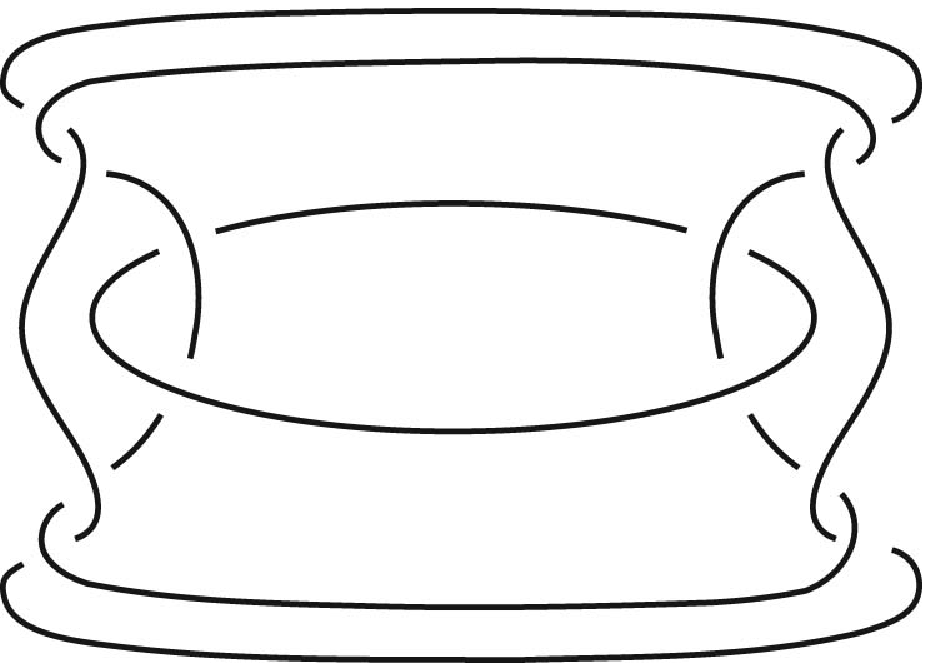}
{\scriptsize
    \put(-290,78){$l'_1$}
    \put(-290,46){$l'_2$}
    \put(-288,14){$l_3$}
    \put(-135,45){$l'_1$}
    \put(-98,45){$l'_2$}
    \put(-70,21){$l_3$}
    \put(-285,-10){R=1/4}
    \put(-65,-12){R=1/8}}
\vspace{.55cm} \caption{}
\end{figure}

This concludes the proof of the relative-slice problem. It remains to show that the slices
$S_1, S_2$ for $l_1,\, l_2$ constructed above are standard. We start by recording the data
involved in their construction. $({\gamma}, {\delta})$ is a Hopf link in $\partial D^4$,
$D'=D\smallsetminus$(collar on $\partial D^4)$. Since the slices were described in terms of the
radial Morse function on $D'$, to be specific consider $D'$ as the $4-$ball of radius $1$ in $D^4$
of radius $2$. The curves $l_1,\, l_2$ are in
the boundary of $D'$; extending them by the product $l_i\times[1,2]$ we will consider
them as curves in $\partial D^4$.

For the rest of this argument, we only need to consider
the curves $l_1,\, l_2$ and their slices; the slices
for the other components are disregarded. $l_1,\, l_2$ form the unlink and therefore bound disjoint
embedded disks $D_1, D_2$ in $\partial D^4$. We will show that the slices $S_1, S_2$ for $l_1,\, l_2$
are standard by constructing disjoint embedded $3-$balls $B_1, B_2$ in $D^4$, with $\partial B_i=D_i
\cup S_i$ for each $i$. The existence of these $3-$balls provides an isotopy in $D^4$ from $S_1, \, S_2$ to
$D_1, \, D_2$ and shows that the slices are standard.

\begin{figure}[ht]
\vspace{.5cm}
\includegraphics [width=3.8cm]{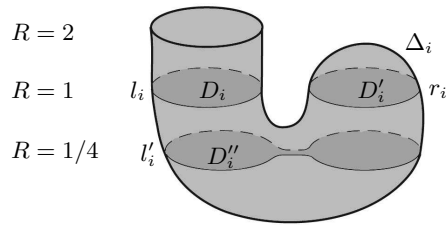}
{\scriptsize
    \put(-160,72){$R=2$}
    \put(-160,50){$R=1$}
    \put(-160,27){$R=1/4$}
    \put(-115,50){$l_i$}
    \put(-2,51){$r_i$}
    \put(-111,26){$l'_i$}
    \put(-89,50){$D_i$}
    \put(-29,50){$D'_i$}
    \put(-86,25){$D''_i$}
    \put(-11,68){${\Delta}_i$}}
\vspace{.45cm} \caption{The $3-$ball $B_i$.}
\end{figure}

\begin{figure}[ht]
\vspace{.5cm}
\includegraphics [width=12cm]{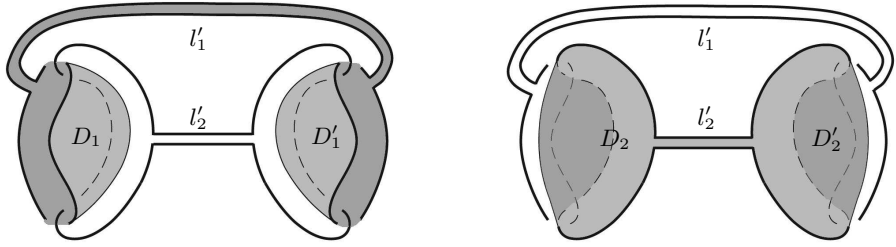}
{\scriptsize
    \put(-270,76){$l'_1$}
    \put(-78,46){$l'_2$}
    \put(-315,39){$D_1$}
    \put(-224,39){$D'_1$}
    \put(-115,38){$D_2$}
    \put(-35,38){$D'_2$}
    \put(-78,76){$l'_1$}
    \put(-270,46){$l'_2$}}
\vspace{.45cm} \caption{Disjoint disks $D''_1, D''_2$: on the left
$D''_1$ is a band sum of $D_1, \, D'_1$, on the right $D''_2$ is a
band sum of $D_2, \, D'_2$.}
\end{figure}

The construction of $S_i$ is illustrated in figure 15. The vertical axis in this figure
corresponds to the radial component in $D^4$. There is a single maximum point given by the
core ${\Delta}_i$ of the $2-$handle attached to $D'$ along $r_i$. Recall that
${\Delta}_1$, ${\Delta}_2$ are embedded in $D^4$ in
a standard way, and so they are isotopic to disjoint embedded disks $D'_1$, $D'_2$ bounded by
$r_1$, $r_2$ in the $3-$sphere
slice $\partial D'=S^3_1$. The curves $l_i,\,  r_i$, and the disks bounded by them: $D_i, \, D'_i$ move
by an isotopy as $R$ decreases from $1$ to $1/4$ until the index $1$ critical points of the slices
at $R=1/4$ (shown in figure 14).

The analysis of the disks at the level $R=1/4$ is presented in figure 16. At the level of these
critical points, the disks $D_i$ and $D'_i$ are band-summed, and the result: the disks $D''_1, \, D''_2$ are disjoint.
The component $l'_1$ on the left in figure 16 bounds $D''_1$, the component $l'_2$ on the right bounds $D''_2$. (Figure 16
has two copies of the link $(l'_1,l'_2)$ just for convenience of visualization of the disks
$D''_1$, $D''_2$.) Finally, at $R<1/4$ the disks $D''_i$ move by an isotopy and shrink to points.

We summarize the construction of the disjoint $3-$balls $B_i$, $i=1,2$: in the $3-$sphere
$S^3_R$, each component of $S^3_R\cap S_i$ bounds a disk: $l_i=\partial D_i$, $r_i=\partial D'_i$,
$l'_i=\partial D''_i$. These disks are the levels of the radial Morse function restricted to $B_i$.
This concludes the proof that $(A,{\alpha})$ embeds into $(S_0,{\gamma}_0)$. \qed

\begin{figure}[ht]
\vspace{.5cm}
\includegraphics[height=3.5cm]{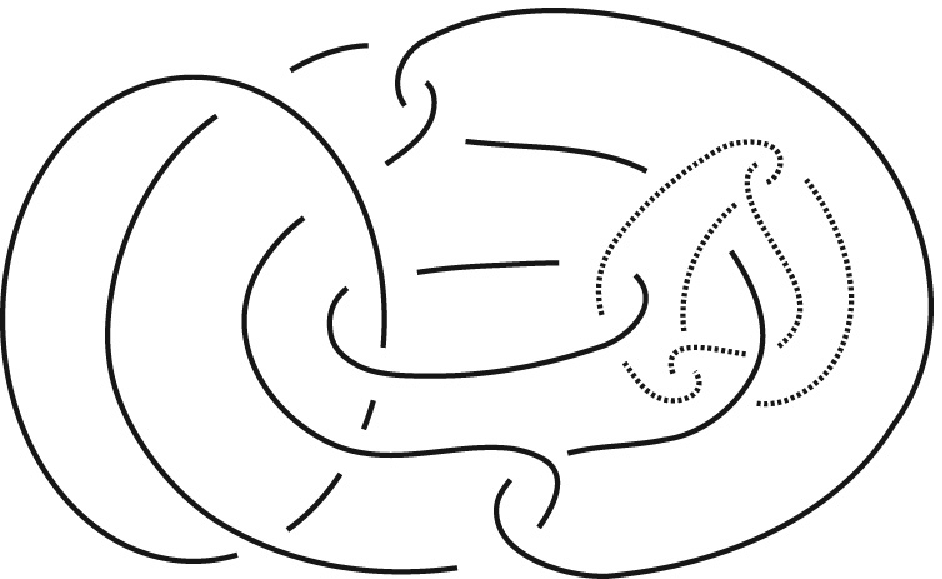} \hspace{1.2cm}  \includegraphics[height=3.5cm]{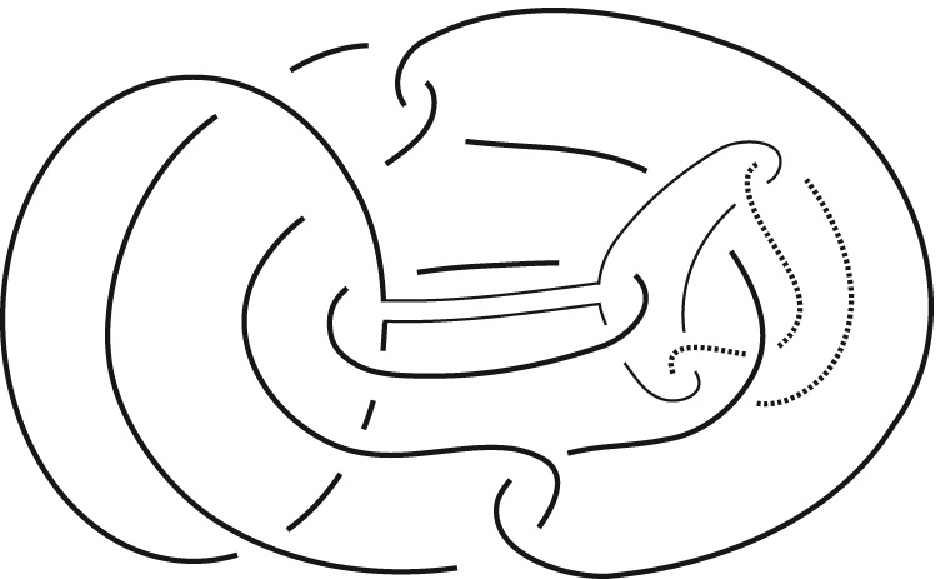}
{\scriptsize
    \put(-372,65){$l_1$}
    \put(-300,-10){$l_2$}
    \put(-228,-2){$l_3$}
    \put(-282,42){$l_4$}
    \put(-220,31){$r_2$}
    \put(-248,78){$r_1$}}
\vspace{.45cm}
\caption{}
\end{figure}

\begin{lemma} \label{lemma B} \sl $B$ embeds in a
collar on its attaching curve. More precisely, there
exists a proper embedding $(B, {\beta})\subset
(S^1\times D^2\times [0,1],\, S^1\times \{0\}\times\{0\})$.
\end{lemma}

\begin{figure}[ht]
\vspace{.5cm}
\includegraphics[height=3.5cm]{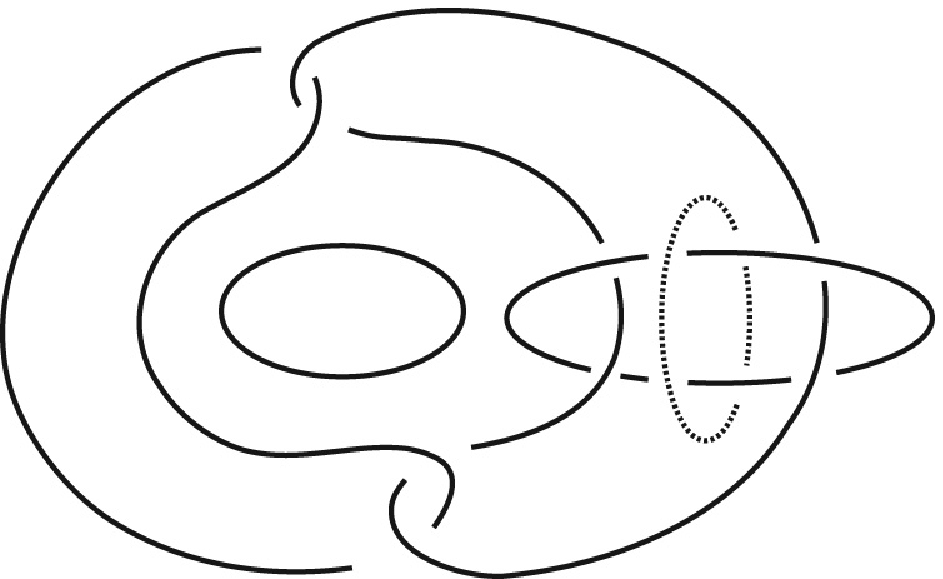} \hspace{1.2cm}  \includegraphics[height=3.5cm]{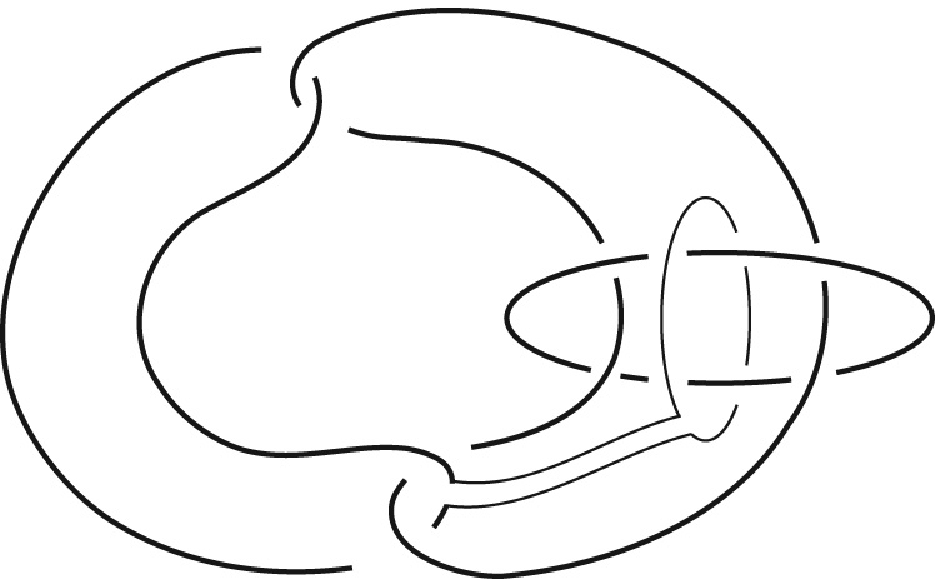}
{\scriptsize
    \put(-211,55){$l_1$}
    \put(-320,-10){$l_2$}
    \put(-245,-2){$l_3$}
    \put(-309,42){$l_4$}
    \put(-248,70){$r_2$}}
\vspace{.45cm} \caption{}
\end{figure}

One needs to show that $(B, {\beta})$ embeds in the complement of a
standard disk bounded by the meridian to ${\beta}$. The proof is
again a relative-slice problem, shown in figure 17. Here $l_1$ is
the meridian which is required to bound a standard disk; $l_2,l_3,
l_4$ are the attaching curves of the $2-$handles of $B$, and $r_1,
r_2$ are the attaching curves for the $2-$handles attached to $D'$.
Therefore the link $l_1,\ldots,l_4$ has to be sliced in
$D'\cup_{r_1,r_2}$(zero-framed $2-$handles), so that the slice for
$l_1$ is standard in $D^4$.

Taking a connected sum of $l_1$ and
$r_1$ as shown in figure 17, one gets the link on the left in figure 18. Now
taking a connected sum of $l_2$ and $r_2$ results in the trivial
link, and the components are capped off with  disjoint disks in
$D'$. The proof that the slice for $l_1$ is standard is directly analogous to the
corresponding proof in lemma \ref{lemma A}. \qed

\begin{figure}[ht]
\vspace{.5cm}
\includegraphics[height=3.5cm]{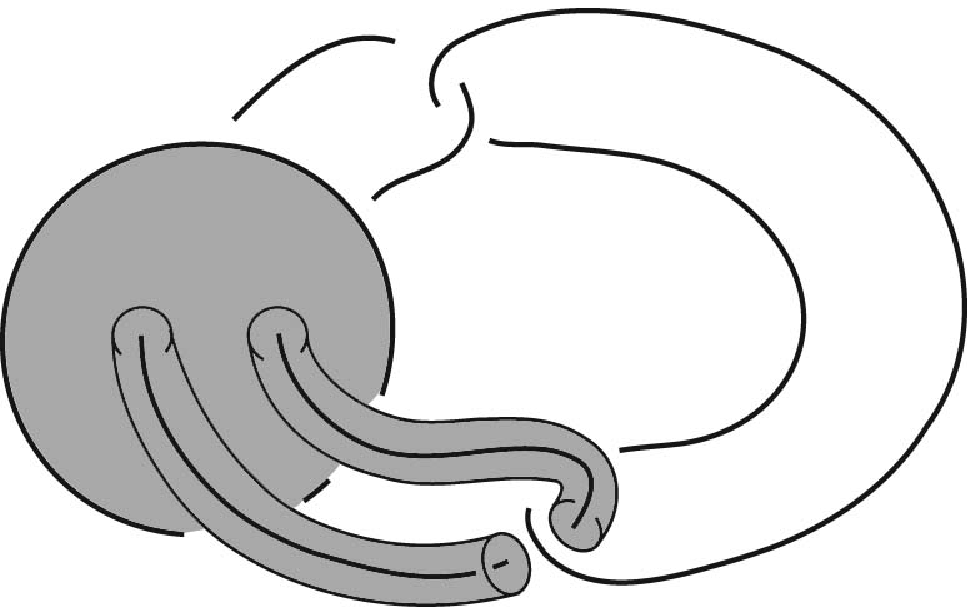}
{\small
    \put(-180,42){$l_1$}
    \put(-126,91){$l_2$}
    \put(-6,78){$l_3$}}
\vspace{.45cm} \caption{}
\end{figure}

{\em Proof of theorem} \ref{theorem} in the central case $L=Bor$,
the Borromean rings, follows from lemmas \ref{lemma A}, \ref{lemma
B}. The components $l_i$ of $Bor$ bound in $D^4$ disjoint embedded
surfaces $S_i$:  $S_1$ is a genus one surface, and $S_2, S_3$ are
disks. Thinking of the radial coordinate of $D^4$ as time where
$\partial D^4$ corresponds to time $0$, $l_1$ bounds a surface $S_1$
(shaded in figure 19) at time $1/2$ and the other two components
bound disjoint disks at $t>1/2$. Consider three decompositions of
$D^4$: $(A_1,B_1)$ is the decomposition constructed in section
\ref{decompositions}. Define $(A_2, B_2)$ and $(A_3, B_3)$ to be the
trivial decomposition: $A_2=A_3=D^2\times D^2$, $B_2=B_3$ are
collars on their attaching curves. Lemmas \ref{lemma A}, \ref{lemma
B} imply that the Borromean rings are weakly $A-B$ slice with these
decompositions.

To prove theorem \ref{theorem} for all links with trivial linking
numbers, a variation of lemmas \ref{lemma A}, \ref{lemma B} is
needed, for higher genus surfaces. That is, given any $g$ there is a
decomposition $D^4=A_g\cup B_g$ such that $A_g$ embeds in (surface
of genus $g)\times D^2$ and $B_g$ embeds in a  collar. These are
variations of the decompositions $A, B$ in definition
\ref{definition}; the case $g=2$ is shown in figure 20. The proof is
analogous to the proof of lemmas \ref{lemma A}, \ref{lemma B}. To
complete the proof of theorem \ref{theorem}, note that the
components of any link with trivial linking numbers bound disjoint
embedded surfaces in $D^4$. \qed

\begin{figure}[ht]
\vspace{.5cm}
\includegraphics[height=5.5cm]{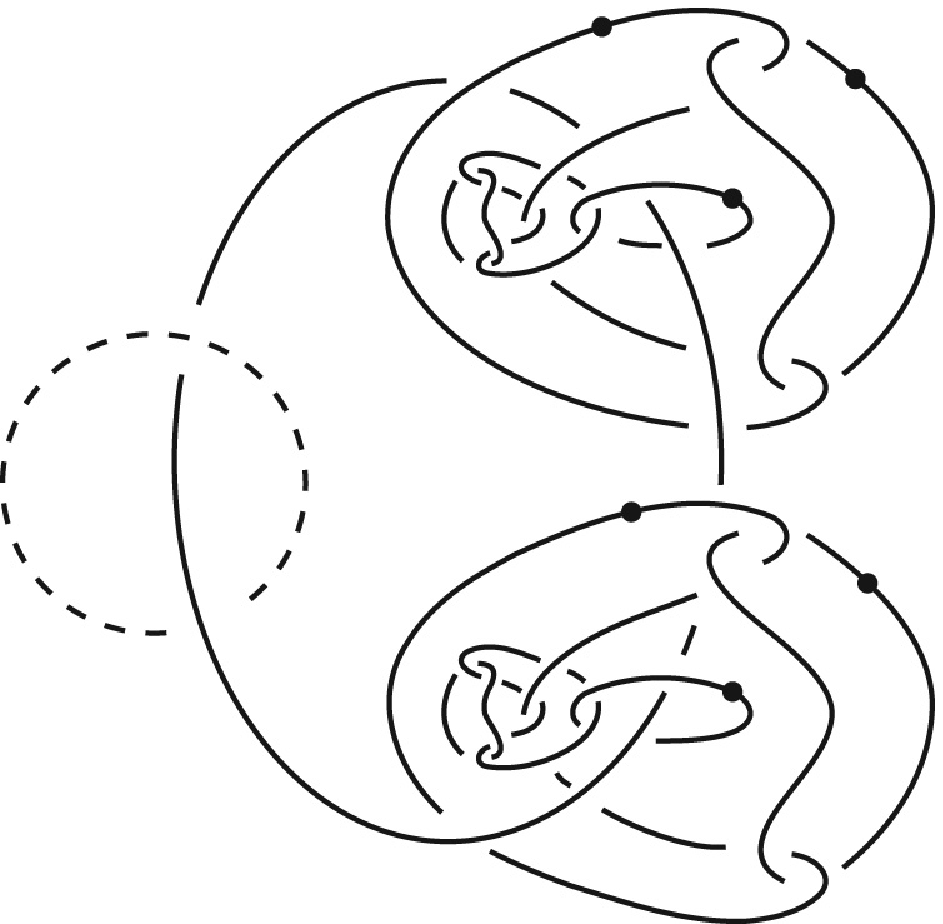} \hspace{2cm} \includegraphics[height=5.5cm]{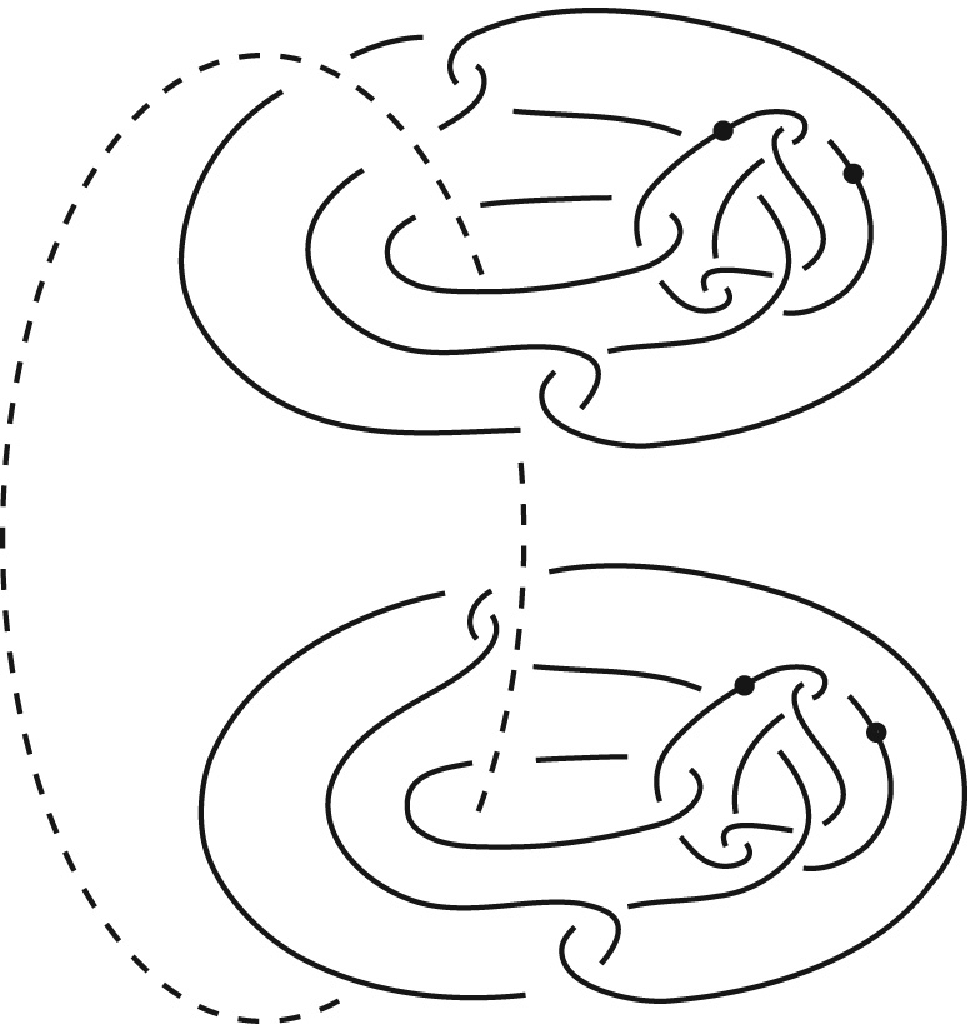}
{\scriptsize
    \put(-318,8){$0$}
    \put(-304,115){$0$}
    \put(-288,101){$0$}
    \put(-304,32){$0$}
    \put(-288,17.5){$0$}
    \put(-81,-7){$0$}
    \put(-41,-7){$0$}
    \put(-63,30.5){$0$}
    \put(-83,80){$0$}
    \put(-43,80){$0$}
    \put(-66,117){$0$}}
    \vspace{.45cm} \caption{}
\end{figure}

{\em Remark.} If the embeddings $(A, {\alpha})\subset (S_0,
{\gamma}_0)\subset (D^4,S^3)$ and $(B, {\beta})\subset (M,
{\gamma})\subset (D^4, S^3)$, constructed in lemmas \ref{lemma A},
\ref{lemma B} were {\em standard}, then taking the complement of the
six submanifolds (three copies of each of $A$ and $B$) bounding the
Borromean rings and their parallel copy in $D^4$ and gluing up the
boundary one would get a solution to a canonical surgery problem.
Considering the generalized Borromean rings, one would get solutions
to all canonical problems, and therefore a proof of the topological
$4-$dimensional surgery conjecture for all fundamental groups.
However the embeddings constructed in the proof are not standard.
This raises the question mentioned in the introduction: Given a
decomposition $D^4=A\cup B$, is one of the {\em embeddings}
$A\hookrightarrow D^4$, $B\hookrightarrow D^4$ necessarily robust?

\bigskip

\section{A decomposition of the $4-$sphere} \label{sphere}

The two parts $A,B$ of any decomposition of $D^4$  have inherently
different properties. For example, due to Alexander duality the
attaching circle ${\alpha}, {\beta}$ on exactly one side vanishes in
its rational first homology group. However it turns out that when
completed to a decomposition of the $4-$sphere, the construction in
the proof of theorem \ref{theorem} is quite symmetric: the
decomposition $D^4=A\cup B$ extends to $S^4=B\cup B$. We record this
observation here since it seems likely that this symmetry plays a
role in the properties of the decomposition used above.

\begin{figure}[ht]
\vspace{.5cm}
\includegraphics[height=3.5cm]{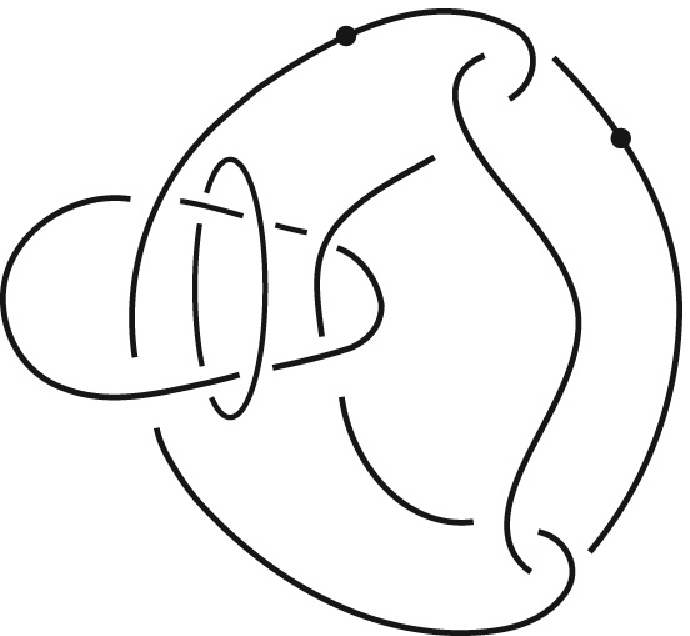} \hspace{2cm}  \includegraphics[height=3.5cm]{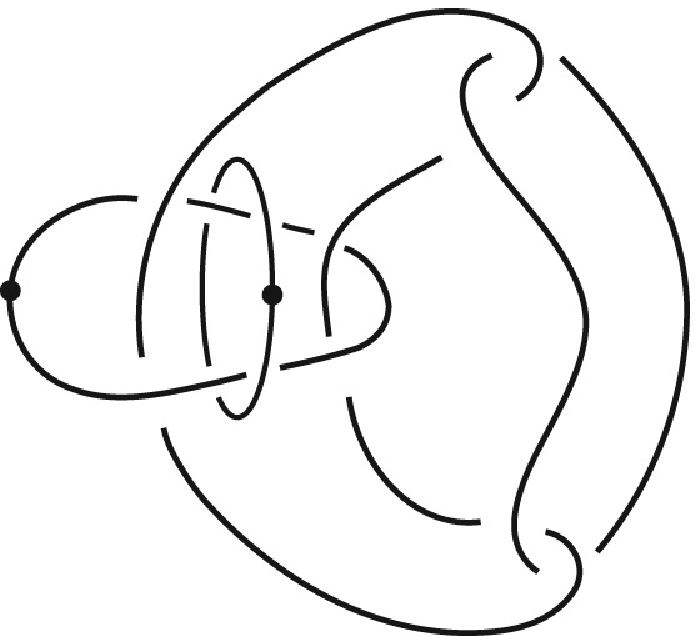}
{\scriptsize
    \put(-243,25){$0$}
    \put(-272,26){$0$}
    \put(-295,12){$A_1$}
    \put(-110,15){$B$}
    \put(-67,93){$0$}
    \put(-9,81){$0$}}
\vspace{.35cm} \caption{}
\end{figure}

Consider the decomposition $D^4=A\cup B$ constructed in definition
\ref{definition}, see figures 7, 8. The attaching curves ${\alpha},
{\beta}$ form the Hopf link in $\partial D^4$. Consider
$S^4=D^4\cup($another copy $D'$ of the $4-$ball$)$, and let
${\alpha}$ bound the standard disk in $D'$. In terms of the
handlebodies, we attach a zero-framed $2-$handle (a thickening of
this disk) to $A$ along ${\alpha}\times D^2$. Denote the result by
$A_1$. Considering the handle diagram for $A$ in figures 8, 12, one
sees that the Kirby diagram for $A_1$ is given in the first part of
figure 21.

\begin{figure}[ht]
\vspace{.5cm}
\includegraphics[height=2.9cm]{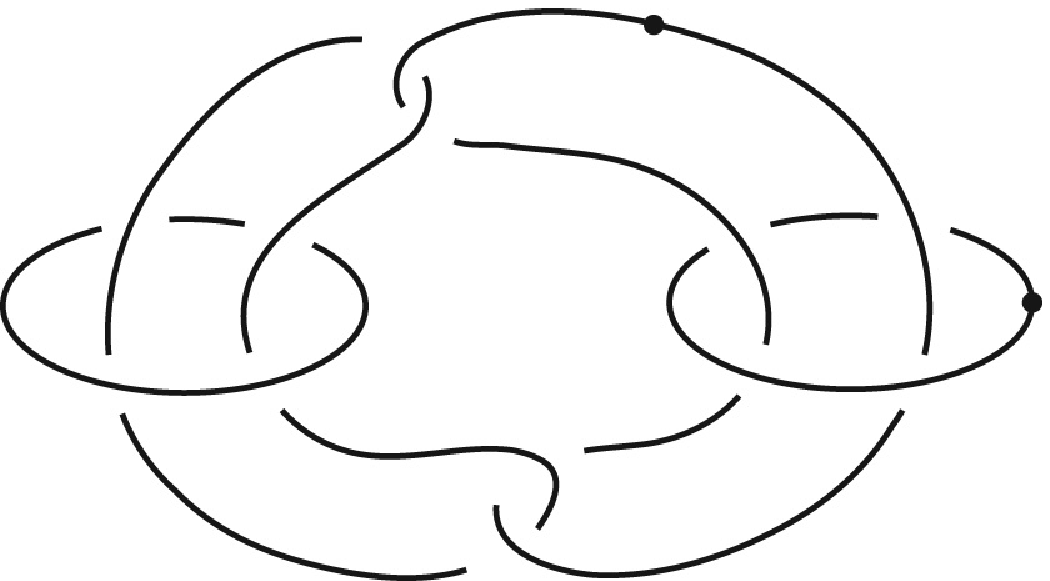}
{\scriptsize
    \put(-125,70){$0$}
    \put(-158,35){$0$}}
    \vspace{.35cm} \caption{}
\end{figure}

Note that the complement of $A_1$ in $S^4$ is homeomorphic to $B$,
since the complement to the disk bounded by ${\alpha}$ in $D'$ is
just a collar on ${\beta}$. The handle diagram for $B$ is presented
in the second part of figure 21. (The dashed circle in figure 7 has
to be replaced by a circle with the dot, since the diagrams in
figure 21 are drawn in $S^3$, not in the solid torus. Then a $1-,
2-$handle pair canceled, and one gets the diagram in figure 21.) The
proof is concluded by the observation that this is a symmetric link:
both diagrams in figure 21 are isotopic to the one in figure 22. \qed

\end{document}